\documentclass{svjour3}
\smartqed  % flush right qed marks, e.g. at end of proof
\usepackage{graphicx}
\usepackage{latexsym,bm,amsmath,amssymb}
\usepackage{rotating}
\usepackage{subcaption}
\usepackage{multirow,booktabs,cases}
\usepackage[misc]{ifsym}
\usepackage[style=1]{mdframed}
\usepackage{algorithm,algorithmic}
\usepackage[colorlinks=true]{hyperref}
\hypersetup{urlcolor=blue, citecolor=blue,linkcolor=blue}
\usepackage{epstopdf}
\usepackage{bm}
\newcommand{\R}{\mathbb{R}}
\def\x{{\bf x}}
\def\y{{\bf y}}
\def\z{{\bf z}}
\def\w{{\bf w}}

\def\c{{\bf c}}
\def\d{{\bf d}}

\def\v{{\bf v}}
\def\t{{\bf t}}
\def\e{{\bf e}}
\def\0{{\bf 0}}
\def\A{{\cal A}}
\def\B{{\cal B}}

\allowdisplaybreaks

\begin{document}
\graphicspath{{./PIC/}}

\title{Biquadratic Cauchy Tensors and Spherical Biquadratic Polynomial Programming}
\titlerunning{Biquadratic Cauchy tensor and Biquadratic Programming}     % if too long for running head

\author{Haibin Chen \and Yixuan Chen\and Liqun Qi}

%\authorrunning{H. Chen \and H. He} % if too long for running head

\institute{Authors are listed alphabetically. \\
H. Chen \at School of Management, Qufu Normal University, Rizhao, Shandong, China, 276800. \\
\email{chenhaibin508@qfnu.edu.cn}
\and Y. Chen \at School of Management, Qufu Normal University, Rizhao, Shandong, China, 276800. \\
\email{yixuanchen1201@163.com}
\and L. Qi \at Department of Applied Mathematics, The Hong Kong Polytechnic University, Hung Hom, Kowloon, Hong Kong, China.\\
\email{maqilq@polyu.edu.hk}
%\and G. Zhou \at Curtin Centre for Optimisation and Decision Science,
%School of Electrical Engineering, Computing and Mathematical Sciences (EECMS),
%Curtin University, Perth, WA 6845, Australia. \\
%\email{g.zhou@curtin.edu.au}
}

\date{Received: date / Accepted: date}
% The correct dates will be entered by the editor

\maketitle

\begin{abstract}
This paper addresses biquadratic polynomial programming (BPP), an NP-hard optimization problem closely related to biquadratic tensors. We first establish several necessary and sufficient conditions for the positive semi-definiteness and positive definiteness of biquadratic Cauchy tensors. Leveraging the structured properties of these tensors, we then prove that the BPP and its equivalent multilinear formulation share the same set of optimal solutions. This result allows us to establish the global sequence convergence of the proximal alternating minimization (PAM) algorithm via the Kurdyka-{\L}ojasiewicz (KL) property, extending the analysis in \cite{Chen2022}. Furthermore, by reformulating the equivalent multilinear problem as an unconstrained optimization model, we enable the analysis of its KL exponent and derive an explicit expression for the convergence rate of PAM. Finally, numerical experiments are conducted on both biquadratic Cauchy tensors and general biquadratic tensor instances to evaluate the efficiency, stability, and practical performance of the proposed algorithm.
% \PACS{PACS code1 \and PACS code2 \and more}
\end{abstract}

%%%%%%%%%%%%%%%%%%%%%%%%%%%%%%%%%%%%%%
\section{Introduction}\label{Int}
Polynomial optimization over spherical constraints is a class of challenging non-convex problems that arises in multiple fields, including signal processing, quantum mechanics, and multivariate statistical analysis \cite{bar2007,CD2007,JB01,QT03,Zhou12}. In this paper, we consider one of a fundamental task, i.e., biquadratic polynomial programming (BPP):
$$
\min f(\x,\y)~~~\mbox{s.t.}~  ||\x||=1, ||\y||=1, \x\in\mathbb{R}^m, \y\in\mathbb{R}^n,
$$
which involves optimizing a homogeneous polynomial of degree four and paired blocks subject to each block lying on a unit sphere.
The BPP frequently arises from the strong ellipticity condition problem
in solid mechanics and the entanglement problem in quantum physics (e.g., see \cite{Han2009,WA1996}) and
the best rank-one approximation of a tensor (e.g., see \cite{Qi2005,QW2009,WCW2018}), which plays an important role
in wireless communication systems, data analysis, and independent component analysis.

On one hand, the biquadratic polynomial $f(\x,\y)$ has a close relationship with biquadratic tensors \cite{CQ2026,QC2025,QCX2025,QCX20252,QHZ2021}. The biquadratic tensor is a kind of structure tensors, which have been systematically studied since its wide applications in statistics, mechanics, the theory of relativity, and polynomial theory \cite{QDH2009,QHZ2021,Zhao23}.
Based on this relationship, researchers in the field of numerical algebra use it to compute the M-eigenvalues of tensors, as well as address scientific problems such as determining the strong ellipticity condition in elasticity theory \cite{Han2009,LLL2019,Li2024,WQZ2009}.
On the other hand, it has been proved that the BPP is NP-hard \cite{Ling09},which implies that finding its global optimal solutions is not an easy task in general cases.
Therefore, researchers are encouraged to develop some effective approximation algorithms for finding approximately optimal solutions, e.g., see \cite{HLZ10,Ling17,Ling09,LZ10,Nie15,So11,ZLQ11} and references therein.

To the best of our knowledge, the semi-definite programming (SDP) is one of
the most popular tools for polynomial optimization problems. For instance, Ling et al. \cite{Ling09}
have proposed various SDP approximation methods to solve BPP with unit spherical
constraints, in addition to analyzing their approximation bounds under a slightly weaker
approximation notion. When dealing with BPP equipped with generic quadratic constraints,
Zhang et al. \cite{ZLQ11} proposed an SDP relaxation to the underlying BPP and extended
their approximation bound results to the cases where decision variables are complex.
Later, Yang and Yang \cite{YY2012} gave a systematical study on BPPs under several different constraints.
Moreover, they presented several variational approaches and gave a provable estimation for
the approximation solutions. For more results on SDP and its applications in polynomial optimization, see \cite{HNY2024,Nie15,Nie15a,NT2024,Nie,Nie2024}.
However, the SDP problems required to be solved by the SOS approach grow exponentially, making the SOS approach only
viable for low dimensional problems.

Unlike the SDP methods, Hu and Huang \cite{HH2011} proposed a
state-of-the-art first order method, i.e., alternating direction method of multipliers (ADMM),
for solving their SDP relaxation model of BPPs. One remarkable property of the ADMM is that it is a globally convergent algorithm for BPPs.
Then, Chen et al. \cite{CSLZ2012} proposed a novel maximum block improvement (MBI) method for directly solving spherically constrained homogeneous polynomial optimization problems. Their method is guaranteed to converge to a stationary point and is applicable to any optimization model with separable block constraints. Empirically, MBI performs stably on polynomial optimization problems when initialized with a reasonably good starting point. By exploiting the multi-block structure of polynomial functions, Jiang et al. \cite{Jiang} reformulated the problem as a linearly constrained model, which facilitates the use of ADMM. As shown in \cite{Jiang}, ADMM is straight forward to implement because each subproblem admits a closed-form solution. However, its convergence analysis relies on a relatively strong condition (namely, $\lim_{k\rightarrow\infty}\|\z^{(k+1)}-\z^{(k)}\|=0$, assumed in \cite{Jiang}), which can be difficult to verify for highly nonlinear and nonconvex polynomial problems. To relax the theoretical requirements, Wang et al. \cite{Wang2015} introduced a new block improvement method (BIM)-essentially a shifted power method [11]-that partitions variables into two blocks. The efficiency of BIM was demonstrated through numerical experiments.
A common feature of the above three methods is their block-updating scheme. As noted in \cite{CSLZ2012}, MBI is inspired by the block coordinate descent (BCD) method \cite{RHL13,Tseng2001}, but differs in that it updates only the block yielding the maximum improvement at each iteration, similar to the Gauss-Southwell rule. In contrast, ADMM updates all blocks sequentially (or alternately) using an augmented Lagrangian framework. Computational results in \cite{Jiang} indicate that ADMM often runs faster than MBI. Meanwhile, BIM \cite{Wang2015}, as an improvement over MBI, updates two block variables alternately with fixed step sizes at each iteration, resembling a gradient-type method. Consequently, it does not guarantee optimality for each subproblem.

Recently, Chen et al. \cite{Chen2022} investigated a class of fourth-degree polynomial optimization problems, where an equivalent augmented multi-linear model incorporating a regularization term is introduced. This augmentation renders the model concave while preserving the equivalence between a multilinear polynomial optimization with the original BPP under unit spherical constraints. Furthermore, a practical method, i.e, proximal alternating minimization (PAM) algorithm was provided. Notably, the PAM algorithm \cite{Chen2022} can be viewed as an application of the BCD method \cite{Tseng2001} with a cyclic rule and proximal terms to polynomial optimization models. In contrast to the MBI method  \cite{CSLZ2012}, the proposed PAM employs a Gauss-Seidel (alternating) rule to update block variables, allowing efficient use of the latest information for potentially better improvement.
Above all, numerical evidence from these studies confirms that block-updating approaches perform well for polynomial optimization problems. Nevertheless, whether the underlying multi-block (i.e., more than two blocks) structure has been fully exploited remains an open question, motivating us to pursue further research in this direction for polynomial models.

In this paper, we further investigates the PAM algorithm for BPPs with spherical constraints. Prior to that, the positive semi-definiteness of biquadratic Cauchy tensors is examined, which serves as a special numerical case for validating the performance of PAM. The main contributions of this paper are summarized as follows:

\noindent{(1)} In \cite{QCX20252}, the conditions for a biquadratic Cauchy tensor to be completely positive were examined. While a biquadratic completely positive tensor is necessarily positive semi-definite (PSD), the converse does not always hold. This paper establishes several necessary and sufficient conditions for the positive semi-definiteness and positive definiteness of biquadratic Cauchy tensors (see Theorems \ref{th1} and \ref{th2}). Notably, all the proposed conditions can be efficiently verified directly from the generating vectors of the tensor. Furthermore, the monotonicity of the biquadratic polynomial associated with such tensors is also investigated (see Theorems \ref{th3} and \ref{th4}).

\noindent{(2)} By the structure property (negative semi-definite) of the biquadratic tensors, the PAM algorithm proposed in \cite{Chen2022} relies primarily on the equivalence between the multilinear formulation and the original BPP in terms of sharing the same optimal value. However, it remains an open question how to recover an optimal solution of the BPP from the optimal solutions of the multilinear problem. Moreover, the global sequence convergence of the algorithm was not established there. In this paper, we strengthen this equivalence by proving that the two problems also share the same optimal solutions. Furthermore, we establish global sequence convergence via the Kurdyka-{\L}ojasiewicz (KL) property. Additionally, the equivalent multilinear problem is reformulated as an unconstrained optimization problem, which enables the analysis of the KL exponent and yields an explicit expression for the convergence rate of PAM.

\noindent{(3)}  Several numerical experiments are conducted to evaluate the practical performance of PAM. Both biquadratic Cauchy tensors and general biquadratic tensor instances are considered to examine the efficiency and numerical stability of the algorithm. Moreover, we compare PAM with MBI and ADMM, which are two representative block-type algorithms widely used for polynomial optimization. Finally, a sensitivity analysis is performed on parameters $\alpha$ and $\gamma$. The results indicate that while a smaller parameter often accelerates convergence, it may also increase the risk of stopping at local solutions. Conversely, a larger parameter tends to improve the likelihood of converging to the minimum eigenvalue, albeit at a higher computational cost.

The remainder of this paper is organized as follows. Section 2 introduces the preliminary knowledge on biquadratic tensors and BPPs. Section 3 presents verifiable necessary and sufficient conditions for their positive (semi-)definiteness. In Section 4, we recall the equivalence between the BPP and the multilinear polynomial optimization. It is
strengthened that the two problems also share the same optimal solutions.
The proposed PAM algorithm and its convergence analysis are provided in Section 5. Section 6 details the experimental setup and reports comparative numerical results. Finally, Section 7 concludes the paper and suggests potential future research directions.

\section{Preliminaries}\label{Sec2}
To make a clear visual distinction, in this paper we use the lowercase letters, boldface letters, and calligraphic uppercase letters to denote scalars, vectors and tensors respectively in general (e.g., the scalar $a$, the decision variable $\x$, and the tensor $\A$). Let $\mathbb{R}^n$ denote the real $n$-dimensional Euclidean space.
For a positive integer $m$, let $[m]=\{1,2,\cdots,m\}$. For $i\in [m], j\in[n]$, let $\e^{(m)}_i$, $\e^{(n)}_j$ be the $i$-th, $j$-th coordinate vector of $\mathbb{R}^m$ and $\mathbb{R}^n$ respectively.

An $m \times n$ biquadratic tensor $\A = (a_{ijkl})$ has real entries $a_{ijkl}$ for $i, k \in [m]$ and $j, l \in [n]$ satisfying
	$$a_{ijkl} = a_{kjil} = a_{ilkj}.$$
Let $\A=(a_{ijkl})$ and $\B=(b_{ijkl})$ be two $m\times n$ biquadratic tensors. The inner product of $\A$ and $\B$ is defined as:
\[
\langle\A,~ \B\rangle=\sum_{i,k\in[m];j,l\in[n]}a_{ijkl}b_{ijkl}.
\]
Accordingly, the induced Frobenius norm is denoted by $\|\A\|_F=\sqrt{\sum_{i,k=1}^m\sum_{j,l=1}^n}a^2_{ijkl}$.
Let $\x \in \mathbb{R}^m$ and $\y \in \mathbb{R}^n$.  Then
	$$f(\x, \y) \equiv \A\x\y\x\y :=\langle\A, \x\circ\y\circ\x\circ\y\rangle= \sum_{i, k=1}^m \sum_{j, l=1}^n a_{ijkl}x_iy_jx_ky_l,$$
where $\x\circ\y\circ\x\circ\y$ is a rank-1 tensor with entries $(\x\circ\y\circ\x\circ\y)_{ijkl}=x_iy_jx_ky_l$.
If $f(\x,\y) \ge 0 (\leq 0)$ for all $\x \in \mathbb{R}^m$, $\y \in \mathbb{R}^n$, then we say that $\A$ is {\bf positive (negative) semidefinite}.  If furthermore, $f(\x,\y) > 0 (<0)$ for all $\x \in \mathbb{R}^m$, {$\y \in \mathbb{R}^n$,} satisfying $\|\x \|_2 = \|\y\|_2 = 1$, then we say that $\A$ is {\bf positive (negative) definite}.

\section{Positive semi-definiteness of biquadtratic Cauchy Tensors}\label{Sec3}

Let $\c = (c_1,\cdots,c_m)^\top\in\mathbb{R}^m$ and $\d = (d_1,\cdots,d_n)^\top\in\mathbb{R}^n$ be generating vectors of biquadratic Cauchy tensor $\A=(a_{ijkl})$ with
entries such that
$$
a_{ijkl}=\frac{1}{c_i + c_k + d_j + d_l},~~ i,k\in[m]; j,l\in[n].
$$
Suppose that $\A$ is well-defined, i.e., $c_i + c_k + d_j + d_l\neq0$ for $i, k\in[m]$ and $j, l\in[n]$.

\begin{theorem}\label{th1} Assume that a biquadratic Cauchy tensor $\A$ has generating vectors $\c = (c_1,\cdots,c_m)^\top\in\mathbb{R}^m$ and $\d = (d_1,\cdots,d_n)^\top\in\mathbb{R}^n$. Then $\A$ is positive semi-definite if and only if $c_i+d_j > 0, i\in[m], j\in[n]$.
\end{theorem}
\proof For the necessary direction, since $\A$ is positive semi-definite, then it holds that
$$
\A\e^{(m)}_i\e^{(n)}_j\e^{(m)}_i\e^{(n)}_j=a_{ijij}=\frac{1}{2(c_i+d_j)}\geq0.
$$
This implies that
$c_i+d_j > 0, i\in[m], j\in[n]$.

For the sufficiency direction, consider the following biquadratic polynomial
$$
\begin{aligned}
\A\x\y\x\y&=\sum_{i,k\in[m];j,l\in[n]}a_{ijkl}x_iy_jx_ky_l=\sum_{i,k\in[m];j,l\in[n]}\frac{x_iy_jx_ky_l}{c_i+c_k+d_j+d_l}\\
&=\sum_{i,k\in[m];j,l\in[n]}\int^1_0t^{c_i+c_k+d_j+d_l-1}x_iy_jx_ky_l{\rm d}t\\
&=\sum_{i,k\in[m];j,l\in[n]}\int^1_0t^{c_i+c_k-\frac{1}{2}}x_ix_kt^{d_j+d_l-\frac{1}{2}}y_jy_l{\rm d}t\\
&=\int^1_0\left(\sum_{i\in[m]}t^{c_i-\frac{1}{4}}x_i\right)^2\left(\sum_{j\in[n]}t^{d_j-\frac{1}{4}}y_j\right)^2{\rm d}t\geq0,
\end{aligned}
$$
where $\x\in\mathbb{R}^m, \y\in\mathbb{R}^n$, and $c_i+d_j > 0$ guarantees that the second equality is well defined. Then the desired results follow.
\qed

Next, we present a sufficient condition for the positive definiteness of biquadratic Cauchy tensors.
\begin{theorem}\label{th2} Assume that a biquadratic Cauchy tensor $\A$ has generating vectors $\c = (c_1,\cdots,c_m)^\top\in\mathbb{R}^m$ and $\d = (d_1,\cdots,d_n)^\top\in\mathbb{R}^n$. Then $\A$ is positive definite if and only if
 $c_i+d_j > 0, i\in[m], j\in[n]$, and further $c_{i_1}\neq c_{i_2}$, $d_{j_1}\neq d_{j_2}$ for all $i_1\neq i_2, i_1,i_2\in[m], j_1\neq j_2, j_1,j_2\in[n]$.
\end{theorem}
\proof For the sufficiency, by Theorem \ref{th1}, we know that $\A$ is positive semi-definite if $c_i+d_j > 0, i\in[m], j\in[n]$. To show its positive definiteness, we prove the conclusion by contradiction. Assume that there are nonzero vectors $\x\in\mathbb{R}^m, \y\in\mathbb{R}^n$ such that $\A\x\y\x\y=0$. Then, it follows that
$$
\begin{aligned}
0&=\A\x\y\x\y=\sum_{i,k\in[m];j,l\in[n]}\int^1_0t^{c_i+c_k+d_j+d_l-1}x_iy_jx_ky_l{\rm d}t\\
&=\int^1_0\left(\sum_{i\in[m]}t^{c_i-\frac{1}{4}}x_i\right)^2\left(\sum_{j\in[n]}t^{d_j-\frac{1}{4}}y_j\right)^2{\rm d}t,
\end{aligned}
$$
which implies that
\begin{equation}\label{e1}
\left(\sum_{i\in[m]}t^{c_i-\frac{1}{4}}x_i\right)\left(\sum_{j\in[n]}t^{d_j-\frac{1}{4}}y_j\right)\equiv0,~t\in(0,1].
\end{equation}
By \eqref{e1} and $c_{i_1}\neq c_{i_2}$, $d_{j_1}\neq d_{j_2}$ for all $i_1,i_2\in[m], j_1,j_2\in[n]$, it holds that
\begin{equation}\label{e2}
\left(x_1+t^{c_2-c_1}x_2+\cdots+t^{c_m-c_1}x_m\right)\left(y_1+t^{d_2-d_1}y_2+\cdots+t^{d_n-d_1}y_n\right)=0,
\end{equation}
for all $t\in(0,1]$. By the continuity of the left function in \eqref{e2}, we know that $x_1y_1=0$.

If $x_1y_1=0$ implies that that $y_1\neq0$. Then we have $x_1=0$ and
$$
\left(x_2+t^{c_3-c_2}x_3+\cdots+t^{c_m-c_2}x_m\right)\left(y_1+t^{d_2-d_1}y_2+\cdots+t^{d_n-d_1}y_n\right)=0,~t\in(0,1].
$$
By the continuity again, it follows that $x_2y_1=0$ and $x_2=0$. Repeat the process above, we obtain that $\x=\0$, and this contradicts the assumption $\x\neq\0$.

If $x_1y_1=0$ implies that $x_1\neq 0$, by a similar analysis to above, we obtain that $\y=\0$, and this contradicts the assumption that $\y\neq\0$.

If $x_1=0,y_1=0$, by \eqref{e2}, we have that
\begin{equation}\label{e3}
\left(x_2+t^{c_3-c_2}x_3+\cdots+t^{c_m-c_2}x_m\right)\left(y_2+t^{d_3-d_2}y_3+\cdots+t^{d_n-d_2}y_n\right)=0,
\end{equation}
for all $t\in(0,1]$. By the continuity of the left function in \eqref{e3}, we know that $x_2y_2=0$. Repeat the proof again, we can still reach a conclusion that contradicts the fact  that $\x$ and $\y$ being non-zero vectors. Therefore, tensor $\A$ is positive definite.

For the necessary condition, we prove the conclusion by contradiction again. Suppose that two elements of $\c$ are equal. Without loss of generality, let $c_1=c_2$. Then take $\x=(1,-1,0,\cdots,0)^\top$ $\in\mathbb{R}^m$ and $\y\in\mathbb{R}^n$. Then, it holds that
$$
\begin{aligned}
\A\x\y\x\y&=\sum_{i,k\in[m];j,l\in[n]}\int^1_0t^{c_i+c_k+d_j+d_l-1}x_iy_jx_ky_l{\rm d}t\\
&=\int^1_0\left(\sum_{i\in[m]}t^{c_i-\frac{1}{4}}x_i\right)^2\left(\sum_{j\in[n]}t^{d_j-\frac{1}{4}}y_j\right)^2{\rm d}t\\
&=\int^1_0\left(t^{c_1-\frac{1}{4}}x_1-t^{c_2-\frac{1}{4}}x_2\right)^2\left(\sum_{j\in[n]}t^{d_j-\frac{1}{4}}y_j\right)^2{\rm d}t\\
&=0,
\end{aligned}
$$
which contradicts the positive definiteness of $\A$, and the desired results hold.
\qed

In the following analysis, we study the monotonicity of biquadratic polynomials corresponding to a biquadratic Cauchy tensor. First of all, for a biquadratic Cauchy tensor
$\A$ with generating vectors $\c\in\mathbb{R}^m$, $\d\in\mathbb{R}^n$, the corresponding biquadratic polynomial is denoted by
$$
f(\x,\y)=\A\x\y\x\y,~\x\in\mathbb{R}^m,\y\in\mathbb{R}^n.
$$
For any $\x\geq\bar{\x}$, $\y\geq\bar{\y}$ (``$\geq$" means element-wise), if $f(\x,\y)\geq f(\bar{\x},\bar{\y})$ ($f(\x,\y)\leq f(\bar{\x},\bar{\y})$), then we say that $f(\x,\y)$ is monotone increasing (monotone decreasing respectively). If $f(\x,\y)> f(\bar{\x},\bar{\y})$ ($f(\x,\y)< f(\bar{\x},\bar{\y})$) when $\x\geq\bar{\x}, \x\neq\bar{\x}$, $\y\geq\bar{\y}, \y\neq\bar{\y}$, then we say that $f(\x,\y)$ is strictly monotone increasing (or strictly monotone decreasing respectively).

\begin{theorem}\label{th3}Assume that a biquadratic Cauchy tensor $\A$ has generating vectors $\c = (c_1,\cdots,c_m)^\top\in\mathbb{R}^m$ and $\d = (d_1,\cdots,d_n)^\top\in\mathbb{R}^n$. Then $\A$ is positive semi-definite if and only if $f(\x,\y)$ is monotone increasing in $\mathbb{R}^m_+\times\mathbb{R}^n_+$.
\end{theorem}
\proof For the sufficiency direction, by the monotone increasing of $f(\x,\y)$, we know that, for any $i\in[m], j\in[n]$,
$$
f(\e^{(m)}_i,\e^{(n)}_j)=\A\e^{(m)}_i\e^{(n)}_j\e^{(m)}_i\e^{(n)}_j=\frac{1}{2(c_i+d_j)}\geq 0.
$$
Combining this with Theorem \ref{th1}, we know that $\A$ is positive semi-definite.

For the necessary condition, for any $\x\geq\bar{\x}$, $\x,\bar{\x}\in\mathbb{R}^m_+$, $\y\geq\bar{\y}$, $\y, \bar{\y}\in\mathbb{R}^n_+$, it follows that
$$
\begin{aligned}
&f(\x,\y)-f(\bar{\x},\bar{\y})=\A\x\y\x\y-\A\bar{\x}\bar{\y}\bar{\x}\bar{\y}\\
&=\sum_{i,k\in[m];j,l\in[n]}\frac{x_iy_jx_ky_l}{c_i+c_k+d_j+d_l}-\sum_{i,k\in[m];j,l\in[n]}\frac{\bar{x}_i\bar{y}_j\bar{x}_k\bar{y}_l}{c_i+c_k+d_j+d_l}\\
&=\sum_{i,k\in[m];j,l\in[n]}\frac{x_iy_jx_ky_l-\bar{x}_i\bar{y}_j\bar{x}_k\bar{y}_l}{c_i+c_k+d_j+d_l}\geq0,
\end{aligned}
$$
where the last inequality is due to the fact that all entries of $\A$ is positive.
\qed

Next, we present a necessary condition for the positive definiteness of biquadratic Cauchy tensor.
\begin{theorem}\label{th4}
Let $\A$ be a biquadratic Cauchy tensor with generating vectors $\c\in\mathbb{R}^m, \d\in\mathbb{R}^n$. Suppose
$\A$ is positive definite. Then the biquadratic polynomial $f(\x,\y)$ is strictly monotone
increasing in $\mathbb{R}^m_+\times\mathbb{R}^n_+$.
\end{theorem}
\proof For any $\x\geq\bar{\x}, \x\neq\bar{\x}$, $\y\geq\bar{\y}, \y\neq\bar{\y}$, without loss of generality, there are indices $i_0\in[m], j_0\in[n]$ such that
$x_{i_0}>\bar{x}_{i_0}, y_{j_0}>\bar{y}_{j_0}$. Combining this with the positive definiteness of $\A$, we know that
$$
\begin{aligned}
&f(\x,\y)-f(\bar{\x},\bar{\y})=\A\x\y\x\y-\A\bar{\x}\bar{\y}\bar{\x}\bar{\y}\\
&=\sum_{i,k\in[m];j,l\in[n]}\frac{x_iy_jx_ky_l}{c_i+c_k+d_j+d_l}-\sum_{i,k\in[m];j,l\in[n]}\frac{\bar{x}_i\bar{y}_j\bar{x}_k\bar{y}_l}{c_i+c_k+d_j+d_l}\\
&=\frac{x_{i_0}^2y_{j_0}^2-\bar{x}_{i_0}^2\bar{y_{j_0}}^2}{2(c_{i_0}+d_{j_0})}+\sum_{i,k\in[m]\backslash\{i_0\};j,l\in[n]\backslash\{j_0\}}\frac{x_iy_jx_ky_l-\bar{x}_i\bar{y}_j\bar{x}_k\bar{y}_l}{c_i+c_k+d_j+d_l}\\
&\geq\frac{x_{i_0}^2y_{j_0}^2-\bar{x}_{i_0}^2\bar{y_{j_0}}^2}{2(c_{i_0}+d_{j_0})}>0,
\end{aligned}
$$
which implies that $f(\x,\y)$ is strictly monotone increasing.
\qed

Next, we present an example to show that the strictly monotone increasing property
for the polynomial $f(\x,\y)$ is only a necessary condition for the positive definiteness
property of $\A$ but not a sufficient condition.
\begin{example} Let $\A$ be a $2\times 3$ biquadratic Cauchy tensor with generating vectors $\c=(1,1)^\top, \d=(2,2,2)^\top$, i.e.,
$$
a_{ijkl}=\frac{1}{6},~\forall~i, k=1,2; j,l=1,2,3,
$$
and the corresponding biquadratic polynomial is
$$
f(\x,\y)=\A\x\y\x\y=\sum_{i,k\in[2];j,l\in[3]}\frac{x_iy_jx_ky_l}{6}.
$$
From a direct computation and by Theorem \ref{th2}, we know that $f(\x,\y)$ is strictly monotone increasing but not positive definite.
\end{example}

\section{Spherical biquadratic programming}\label{sec04}

In this section, we consider the following spherical biquadratic polynomial programming (BPP):
\begin{equation}\label{e4}
\min_{\x,\y\in\mathbb{D}_1} f(\x,\y)=\min_{\x,\y\in\mathbb{D}_1} \A\x\y\x\y,
\end{equation}
where $\A$ is an $m\times n$ biquadratic tensor, and $\mathbb{D}_1:=\{(\x,\y)\in\mathbb{R}^m\times\mathbb{R}^n ~|~ \|\x\|=1, \|\y\|=1\}$.
Similarly, by the same biquadratic tensor $\A$, we have the following biquadratic multi-linear programming (BMP):
\begin{equation}\label{e5}
\min_{{\bf u},\v,\w,\z\in\mathbb{D}_2} F({\bf u},\v,\w,\z)=\min_{{\bf u},\v,\w,\z\in\mathbb{D}_2} \A{\bf u}\v\w\z,
\end{equation}
where $\mathbb{D}_2:=\{({\bf u},\v,\w,\z)\in\mathbb{R}^m\times\mathbb{R}^n\times\mathbb{R}^m\times\mathbb{R}^n ~|~ \|{\bf u}\|=\|\w\|=1, \|\v\|=\|\z\|=1\}$.

As the equivalence between \eqref{e4} and \eqref{e5} has been established in \cite{Chen2022}, the proofs of Theorems \ref{th5} and \ref{th6} are omitted here.
\begin{theorem}\label{th5}
Let $\A$ be an $m\times n$ biquadratic tensor. If $\A$ is negative semi-definite, then \eqref{e4} and \eqref{e5} are equivalent in the sense that
$$
\min_{\x,\y\in\mathbb{D}_1} f(\x,\y)=\min_{{\bf u},\v,\w,\z\in\mathbb{D}_2} F({\bf u},\v,\w,\z).
$$
\end{theorem}
%\proof First of all, it is easy to know that
%\begin{equation}\label{e6}
%\min_{\x,\y\in\mathbb{D}_1} f(\x,\y)\geq\min_{{\bf u},\v,\w,\z\in\mathbb{D}_2} F({\bf u},\v,\w,\z).
%\end{equation}

%On the other hand, by the negative semi-definiteness of $\A$, we know that, for any $\x\in\mathbb{R}^m, \y\in\mathbb{R}^n$,
%\begin{equation}\label{e7}
%f(\x,\y)=\A\x\y\x\y\leq 0.
%\end{equation}
%By \eqref{e7}, for any ${\bf u},\w\in\mathbb{R}^m$ and $\v,\z\in\mathbb{R}^n$, it holds that
%$$
%\begin{aligned}
%0&\geq f({\bf u}+\w, \v-\z)+f({\bf u}-\w,\v+\z)\\
%&=\A({\bf u}+\w)(\v-\z)({\bf u}+\w)(\v-\z)+\A({\bf u}-\w)(\v+\z)({\bf u}-\w)(\v+\z)\\
%&=2(\A{\bf u}\z{\bf u}\z+\A{\bf u}\v{\bf u}\v+\A\w\z\w\z++\A\w\v\w\v)-8\A{\bf u}\v\w\z\\
%&=2(f({\bf u},\z)+f({\bf u},\v)+f(\w,\z)+f(\w,\v))-8F({\bf u},\v,\w,\z),
%\end{aligned}
%$$
%which means that
%$$F({\bf u},\v,\w,\z)\geq\min\{f({\bf u},\z),f({\bf u},\v),f(\w,\z),f(\w,\v)\}.$$
%By the arbitrariness of ${\bf u},\w\in\mathbb{R}^m$, $\v,\z\in\mathbb{R}^n$, we know that
%$$
%\min_{{\bf u},\v,\w,\z\in\mathbb{D}_2} F({\bf u},\v,\w,\z)\geq\min_{\x,\y\in\mathbb{D}_1} f(\x,\y).
%$$
%Combining this with \eqref{e7}, it confirms the equivalence of problems \eqref{e4} and \eqref{e5}, and the desired results hold
%\qed

Generally speaking, the given tensor may not be negative semi-definite. Therefore, a common way is to subtract a regularization term from the original objective function. We introduce the following BPP:
\begin{equation}\label{e8}
\min_{\x,\y\in\mathbb{D}_1} f_{\alpha}(\x,\y)=\A\x\y\x\y-\alpha\|\x\|^2\|\y\|^2.
\end{equation}
By \eqref{e8}, the corresponding BMP is as below:
\begin{equation}\label{e9}
\min_{{\bf u},\v,\w,\z\in\mathbb{D}_2} F_{\alpha}({\bf u},\v,\w,\z)=\A{\bf u}\v\w\z-\alpha\langle{\bf u},\w\rangle\langle\v,\z\rangle.
\end{equation}
If the constant $\alpha$ is assigned an appropriate value, it can be ensured that the function $f_{\alpha}(\x,\y)$ is less than or equal to zero, i.e., the negative semi-definiteness of tensor $\A$. On the other hand, by definitions of $\mathbb{D}_1, \mathbb{D}_2$, we know that \eqref{e4} and \eqref{e8} share the same optimal solutions. Then, we have the following results.
\begin{theorem}\label{th6}
For an $m\times n$ biquadratic tensor $\A=(a_{ijkl})$, if $\alpha\geq\|\A\|$ in \eqref{e8}, then $\A$ is negative semi-definite. Furthermore, \eqref{e8} and \eqref{e9} are equivalent, i.e.,
$$
\min_{\x,\y\in\mathbb{D}_1} f_{\alpha}(\x,\y)=\min_{{\bf u},\v,\w,\z\in\mathbb{D}_2} F_{\alpha}({\bf u},\v,\w,\z).
$$
\end{theorem}
%\proof By conditions, for any $\x\in\mathbb{R}^m, \y\in\mathbb{R}^n$, it holds that
%$$
%\begin{aligned}
%f_{\alpha}(\x,\y)&=\A\x\y\x\y-\alpha\|\x\|^2\|\y\|^2\\
%&=\left(\sum\limits_{\substack{i,k\in[m];\,j,l\in[n]}}a_{ijkl}x_iy_jx_ky_l\right)-\alpha\|\x\|^2\|\y\|^2\\
%&\leq \sqrt{\sum\limits_{i,k\in[m];\, j,l\in[n]}a^2_{ijkl}\sum\limits_{i,k\in[m];\, j,l\in[n]}x^2_iy^2_jx^2_kz^2_l}-\alpha\|\x\|^2\|\y\|^2\\
%&=\|\A\|_F\|\x\|^2\|\y\|^2-\alpha\|\x\|^2\|\y\|^2\\
%&=\left(\|\A\|_F-\alpha\right)\|\x\|^2\|\y\|^2\leq 0,
%\end{aligned}
%$$
%which implies that $\A$ is negative semi-definite. Combining this with Theorem \ref{th5}, it follows that \eqref{e8} and \eqref{e9} are equivalent, and the desired results hold.
%\qed

From Theorem \ref{th6}, we have the following result.
\begin{corollary}\label{corol1} Let $\alpha\geq\|\A\|_F$. Then, problems \eqref{e4} (or \eqref{e8}) and \eqref{e9} share at least a pair of common optimal solutions.
\end{corollary}
\proof For problem \eqref{e4}, it has at least one optimal solution in $D_1$ since the continuous objective function must have optimal solutions on a compact set.
Suppose that $\x^*, \y^*$ are optimal solution of \eqref{e4}. Then, ${\bf u}=\w=\x^*, \w=\z=\y^*$ are optimal solutions of \eqref{e9} and the desired result holds.
\qed

From Corollary \ref{corol1}, a question can be raised naturally, i.e., how can we recover the optimal solution of \eqref{e4} (or \eqref{e8}) from the optimal solutions of \eqref{e9}? To the best of our knowledge, the question has not been considered before.
\begin{theorem}\label{them7} Let $\alpha\geq\|\A\|_F$. Suppose that ${\bf u}^*, \v^*, \w^*, \z^*\in D_2$ are optimal solutions of \eqref{e9}. Then $({\bf u}^*, \z^*), ({\bf u}^*, \v^*), (\w^*, \v^*), (\w^*, \z^*)\in D_1$ are all optimal solutions of \eqref{e4}, i.e.,
$$
f({\bf u}^*, \z^*)=f({\bf u}^*, \v^*)=f(\w^*, \v^*)=f(\w^*, \z^*)=\min_{\x,\y\in\mathbb{D}_1} f(\x,\y).
$$
\end{theorem}
\proof By conditions, suppose that ${\bf u}^*, \v^*, \w^*, \z^*\in D_2$ are optimal solutions of \eqref{e9}, i.e,
$$
({\bf u}^*, \v^*, \w^*, \z^*)=\arg\min\{ F_{\alpha}({\bf u},\v,\w,\z)~|~{\bf u},\v,\w,\z\in\mathbb{D}_2\}.
$$
From Theorem \ref{th6}, we know that the corresponding coefficient tensor of $f_{\alpha}(\x,\y)$ is negative semi-definite, and it follows that
$$
\begin{aligned}
0&\geq f_{\alpha}({\bf u}^*+\w^*, \v^*-\z^*)+f_{\alpha}({\bf u}^*-\w^*,\v^*+\z^*)\\
&=2(f_{\alpha}({\bf u}^*,\z^*)+f_{\alpha}({\bf u}^*,\v^*)+f_{\alpha}(\w^*,\z^*)+f_{\alpha}(\w^*,\v^*))-8F_{\alpha}({\bf u}^*,\v^*,\w^*,\z^*),
\end{aligned}
$$
which means that
$$4F_{\alpha}({\bf u}^*,\v^*,\w^*,\z^*)\geq f_{\alpha}({\bf u}^*,\z^*)+f_{\alpha}({\bf u}^*,\v^*)+f_{\alpha}(\w^*,\z^*)+f_{\alpha}(\w^*,\v^*).$$
Combining this with \eqref{e8}-\eqref{e9} and the fact that
$$
\min\{f_{\alpha}({\bf u}^*,\z^*),f_{\alpha}({\bf u}^*,\v^*),f_{\alpha}(\w^*,\z^*),f_{\alpha}(\w^*,\v^*)\}\geq F_{\alpha}({\bf u}^*,\v^*,\w^*,\z^*),
$$
we obtain that
$$
f_{\alpha}({\bf u}^*,\z^*)=f_{\alpha}({\bf u}^*,\v^*)=f_{\alpha}(\w^*,\z^*)=f_{\alpha}(\w^*,\v^*)=F_{\alpha}({\bf u}^*,\v^*,\w^*,\z^*),
$$
which implies that  $({\bf u}^*, \z^*), ({\bf u}^*, \v^*), (\w^*, \v^*), (\w^*, \z^*)\in D_1$ are all optimal solutions of \eqref{e4}. Therefore the desired result holds.
\qed

\section{Proximal alternating minimization algorithm}\label{sec05}

By Theorem \ref{th6} and the analysis above, we know that finding the optimal value of BPP \eqref{e8} amounts to computing the optimal value of the BMP \eqref{e9}. In this section, due to the multi-block structure of the BMP, we employ the state-of-the-art block coordinate descent method to update all blocks one-by-one in a sequential (a.k.a., Gauss-Seidel) order. Hereafter, we present Algorithm \ref{alg1}, which is called {\it proximal alternating minimization} (PAM) algorithm, for solving problem \eqref{e9}.

\begin{algorithm}[!htbp]
	\caption{(Proximal Alternating Minimization Algorithm for \eqref{e9}).}\label{alg1}
	\begin{algorithmic}[1]
		\STATE Let $\alpha\geq\|\A\|_F$, $\varepsilon>0$, ${\bf u}^{(0)},\v^{(0)},\w^{(0)},\z^{(0)}\in \mathbb{D}_2$, $\gamma_i\geq 0,\;(i=1,2,3,4)$.
		\FOR{$k=0,1,2,\cdots,n$}
		\STATE Update $({\bf u}^{(k+1)},\v^{(k+1)},\w^{(k+1)},\z^{(k+1)})$ sequentially via
		\begin{subnumcases}{\label{PAMA}}
		{\bf u}^{(k+1)}=\arg\min\limits_{\|{\bf u}\|=1}\left[F_{\alpha}({\bf u},\v^{(k)},\w^{(k)},\z^{(k)})+{\frac{\gamma_1}{2}}\|{\bf u}-{\bf u}^{(k)}\|^2\right].\label{PAMA-1} \\
		\v^{(k+1)}=\arg\min\limits_{\|\v\|=1}\left[F_{\alpha}({\bf u}^{(k+1)},\v,\w^{(k)},\z^{(k)})+{\frac{\gamma_2}{2}}\|\v-\v^{(k)}\|^2\right]. \label{PAMA-2}\\
		\w^{(k+1)}=\arg\min\limits_{\|\w\|=1}\left[F_{\alpha}({\bf u}^{(k+1)},\v^{(k+1)},\w,\z^{(k)})+{\frac{\gamma_3}{2}}\|\w-\w^{(k)}\|^2\right]. \label{PAMA-3}\\
		\z^{(k+1)}=\arg\min\limits_{\|\z\|=1}\left[F_{\alpha}({\bf u}^{(k+1)},\v^{(k+1)},\w^{(k+1)},\z)+{\frac{\gamma_4}{2}}\|\z-\z^{(k)}\|^2\right]. \label{PAMA-4}
		\end{subnumcases}
		\STATE Find a pair of $(\x^{(k+1)}, \y^{(k+1)})$ via
		\begin{align*}
		(\x^{(k+1)}, \y^{(k+1)})=\arg\min &\left\{f_{\alpha}({\bf u}^{(k+1)},\v^{(k+1)}),f_{\alpha}({\bf u}^{(k+1)},\z^{(k+1)}), f_{\alpha}(\w^{(k+1)},\v^{(k+1)}),\right. \\
		&\;\;\;\;\left. f_{\alpha}(\w^{(k+1)},\z^{(k+1)})\right\}.
		\end{align*}
		\STATE if {$| f_{\alpha}(\x^{(k+1)},\y^{(k+1)})-f_{\alpha}(\x^{(k)},\y^{(k)})|\leq \varepsilon$}, stop and return a pair of approximate solutions.
		\ENDFOR
	\end{algorithmic}
\end{algorithm}
\begin{remark}\label{r1}
Theoretically, it is possible to set all \(\gamma_i = 0\), as this does not affect the monotonicity of the sequence. However, in practice, we found that choosing \(\gamma_i \in [1,5]\) significantly improves convergence speed and reduces the number of iterations. It is also worth noting that while the convergence proof requires \(\alpha \geq \|\A\|_F\) as a sufficient condition, this is not strictly necessary in practice. To assess the practical impact of these parameters, we conducted additional numerical experiments evaluating the sensitivity of the computed eigenvalue $\lambda$, CPU time, and iteration count with respect to changes in $\gamma_i$. Furthermore, we observed that when $\gamma_i$ is too small, the algorithm often converges to a local minimum (typically the second-smallest tensor eigenvalue). Conversely, larger values of $\gamma_i$ tend to enable convergence to the global minimum eigenvalue but at the cost of increased CPU time. These findings highlight the importance of balancing iteration efficiency and the global convergence property. In particular, we found that setting $\gamma_i$ (or $\alpha$) to be of the same order of magnitude as the Frobenius norm of the tensor yields the robust performance. This observation supports the use of heuristically bounded parameter ranges and may inform future work on adaptive or automated parameter tuning strategies.
\end{remark}

\begin{remark}
(1) Algorithm \ref{alg1} is straightforward to implement under spherical constraints. Since the linearly independent constraint qualification (LICQ) automatically holds, the (local) optimal solution for each subproblem of Algorithm \ref{alg1} is guaranteed to be a KKT point. Considering the subproblem \eqref{PAMA-1} of Algorithm \ref{alg1},
\begin{eqnarray}
\label{e11}
{\bf u}^{(k+1)} = \arg\min_{\|{\bf u}\|=1} F_{\alpha}({\bf u}, \v^{(k)}, \w^{(k)}, \z^{(k)}) + \frac{\gamma_1}{2} \|{\bf u} - {\bf u}^{(k)}\|^2.
\end{eqnarray}
By applying the first-order optimality condition, this subproblem has closed form solution, yielding only two candidates
\begin{eqnarray}
\label{e12}
{\bf u}^{(k+1)} = \pm \frac{\A\v^{(k)}\w^{(k)}\z^{(k)} - \alpha \v^{(k)}\langle\w^{(k)},\z^{(k)}\rangle - \gamma_1 {\bf u}^{(k)}}{\|\A\v^{(k)}\w^{(k)}\z^{(k)} - \alpha \v^{(k)}\langle\w^{(k)},\z^{(k)}\rangle - \gamma_1 {\bf u}^{(k)}\|}.
\end{eqnarray}
Note that if the denominators of \eqref{e12} equal zero, we may set \({\bf u}^{k+1} = {\bf u}^{k}\) and adjust the parameters \(\gamma_1\) and \(\alpha\). Comments regarding the appropriate size of \(\gamma_1\) and \(\alpha\) are provided in Remark~\ref{r1}. Similarly, the KKT points for the subproblems with blocks $\v, \w, \z,$ can be derived analogously. Since all subproblems in Algorithm \ref{alg1} admit finite analytic solutions, the algorithm is computationally efficient for spherically constrained problems.

\noindent(2) Similar algorithms have been applied to various polynomial optimization problems in the literature \cite{CSLZ2012,Chen2022,Wang2015}. For instance, \cite{Chen2022} employs such an algorithm for a biquadratic optimization problem defined on two unit spheres, where only subsequential convergence is established. In contrast, the Maximum Block Improvement (MBI) algorithm \cite{CSLZ2012} updates one block per iteration, while the Block Improvement Method (BIM) \cite{Wang2015} updates two blocks in sequential order. Our proposed algorithm differs by updating all blocks sequentially, which often results in greater improvement in the objective function value for multilinear optimization problems, as demonstrated by computational results in Section 5.

\noindent(3) Furthermore, Step 3 of Algorithm \ref{alg1} can be viewed as a special case of the block coordinate descent (BCD) algorithm \cite{Tseng2001} with a cyclic update rule. For the proximal parameters $\gamma_i (i=1,\ldots,d)$, they can theoretically take any positive values. In our experiments in Section 5, we select $\gamma_i \in (0, 10]$.
\end{remark}

Now, we first establish the subsequential convergence of the sequence $\{\t^{(k)}=({\bf u}^{(k)}, \v^{(k)}, \w^{(k)}, \z^{(k)})\}$ generated by Algorithm \ref{alg1}.
For the sake of simplicity, denote $\t=({\bf u},\v,\w,\z)$ and $F_{\alpha}({\bf u}, \v, \w, \z)=F_{\alpha}(\t)$.
\begin{theorem}\label{thm7}{\rm (Subsequential Convergence)} Let $\{\t^{(k)}\}$ be an infinite sequence generated by Algorithm \ref{alg1}.
Let $\bar{\gamma}=\min\{\gamma_1,\gamma_2,\gamma_3,\gamma_4\}$. Then the following statements hold.

{\rm(i)} For all $k$, the function sequence $\{F_{\alpha}(\t^{(k)})\}$ is nonincreasing and convergent. The sequence $\{\t^{(k)}\}$ satisfies
$\sum_{k=1}^{+\infty}\|\t^{(k+1)}-\t^{(k)}\|<+\infty.$

{\rm(ii)} Suppose $\bar{\t}$ is a cluster point of $\{\t^{(k)}\}$, then $\bar{\t}$ is a KKT point of (\ref{e9}) and
$\lim_{k\rightarrow+\infty}F_{\alpha}(\t^{(k)})=F_{\alpha}(\bar{\t}).$
\end{theorem}
\proof
(i) By Algorithm \ref{alg1}, it follows that
$$
\begin{aligned}
&F_{\alpha}(\t^{(k+1)})+\frac{\gamma_4}{2}\|\z^{(k+1)}-\z^{(k)}\|^2\leq F_{\alpha}({\bf u}^{(k+1)},\v^{(k+1)},\w^{(k+1)},\z^{(k)}),\\
&F_{\alpha}({\bf u}^{(k+1)},\v^{(k+1)},\w^{(k+1)},\z^{(k)})+\frac{\gamma_3}{2}\|\w^{(k+1)}-\w^{(k)}\|^2\leq F_{\alpha}({\bf u}^{(k+1)},\v^{(k+1)},\w^{(k)},\z^{(k)}),\\
&F_{\alpha}({\bf u}^{(k+1)},\v^{(k+1)},\w^{(k)},\z^{(k)})+\frac{\gamma_2}{2}\|\v^{(k+1)}-\v^{(k)}\|^2\leq F_{\alpha}({\bf u}^{(k+1)},\v^{(k)},\w^{(k)},\z^{(k)}),\\
&F_{\alpha}({\bf u}^{(k+1)},\v^{(k)},\w^{(k)},\z^{(k)})+\frac{\gamma_1}{2}\|{\bf u}^{(k+1)}-{\bf u}^{(k)}\|^2\leq F_{\alpha}(\t^{(k)}),
\end{aligned}
$$
which implies
\begin{equation}\label{e13}
F_{\alpha}(\t^{(k+1)})+\frac{\bar{\gamma}}{2}\|\t^{(k+1)}-\t^{(k)}\|^2\leq F_{\alpha}(\t^{(k)}).
\end{equation}
Therefore, the sequence $\{F_{\alpha}(\t^{(k)})\}$ is nonincreasing. Since $F_{\alpha}(\t)$ is bounded on the compact set $\mathbb{D}_2$, we know that $\{F_{\alpha}(\t^{(k)})\}$ is convergent, and the sequence $\{\t^{(k)}\}$ satisfies $\sum_{k=1}^{+\infty}\|\t^{(k+1)}-\t^{(k)}\|<+\infty$ from (\ref{e13}). The conclusion (i) follows.

\noindent(ii) For the sequence $\{\t^{(k)}\}$, it has cluster points by bounded conditions. Without loss of generality, suppose that $\{\t^{(k_j)}\}$ is a subsequence of  $\{\t^{(k)}\}$ with cluster point $\bar{\t}=(\bar{{\bf u}}, \bar{\v},\bar{\w},\bar{\z})$. By the proof of (i), it's clear that $\lim_{k\rightarrow+\infty}F_{\alpha}(\t^{(k)})=F_{\alpha}(\bar{\t})$. On the other hand, from (i) again, it holds that
$$
\lim_{k\rightarrow+\infty}\|\t^{(k+1)}-\t^{(k)}\|=0.
$$
Then, we have that $\lim_{j\rightarrow+\infty}\t^{(k_j+1)}=\bar{\t}.$ Considering the KKT system of the subproblems in Algorithm \ref{alg1}, there are lagrange multipliers $\lambda_1,\lambda_2, \lambda_3, \lambda_4\in\mathbb{R}$ satisfying
$$
\left\{
\begin{aligned}
&\nabla_{{\bf u}} F_{\alpha}({\bf u}^{(k_j+1)},\v^{(k_j)},\w^{(k_j)},\z^{(k_j)})+\gamma_1({\bf u}^{(k_j+1)}-{\bf u}^{(k_j)})-\lambda_1{\bf u}^{(k_j+1)}=\0,\\
&\nabla_{\v} F_{\alpha}({\bf u}^{(k_j+1)},\v^{(k_j+1)},\w^{(k_j)},\z^{(k_j)})+\gamma_2(\v^{(k_j+1)}-\v^{(k_j)})-\lambda_2\v^{(k_j+1)}=\0,\\
&\nabla_{\w} F_{\alpha}({\bf u}^{(k_j+1)},\v^{(k_j+1)},\w^{(k_j+1)},\z^{(k_j)})+\gamma_3(\w^{(k_j+1)}-\w^{(k_j)})-\lambda_3\w^{(k_j+1)}=\0,\\
&\nabla_{\z} F_{\alpha}({\bf u}^{(k_j+1)},\v^{(k_j+1)},\w^{(k_j+1)},\z^{(k_j+1)})+\gamma_4(\z^{(k_j+1)}-\z^{(k_j)})-\lambda_4\z^{(k_j+1)}=\0.\\
\end{aligned}
\right.
$$
Let $j\rightarrow+\infty$, and by the continuity of $F_{\alpha}(\t)$, it follows that
\begin{equation}\label{e14}
\left\{
\begin{aligned}
&\nabla_{{\bf u}} F_{\alpha}(\bar{{\bf u}},\bar{\v},\bar{\w},\bar{\z})-\lambda_1\bar{{\bf u}}=\0,\\
&\nabla_{\v} F_{\alpha}(\bar{{\bf u}},\bar{\v},\bar{\w},\bar{\z})-\lambda_2\bar{\v}=\0, \\
&\nabla_{\w} F_{\alpha}(\bar{{\bf u}},\bar{\v},\bar{\w},\bar{\z})-\lambda_3\bar{\w}=\0, \\
&\nabla_{\z} F_{\alpha}(\bar{{\bf u}},\bar{\v},\bar{\w},\bar{\z})-\lambda_4\bar{\z}=\0.\\
\end{aligned}
\right.
\end{equation}
By the fact that $\bar{{\bf u}},\bar{\v},\bar{\w},\bar{\z}\in\mathbb{D}_2$, we obtain that
$$
\lambda_1=\lambda_2=\lambda_3=\lambda_4=F_{\alpha}(\bar{{\bf u}},\bar{\v},\bar{\w},\bar{\z})=F_{\alpha}(\bar{\t}).
$$
Combining this with (\ref{e14}), it holds that $\nabla F_{\alpha}(\bar{\t})-F_{\alpha}(\bar{\t})\bar{\t}=\0,$
and the desired result holds.
\qed

To prove the global sequence convergence, numerous classical results from the literature have been employed to guarantee it for descent algorithms, including proximal algorithms, forward-backward splitting algorithms, regularized Gauss-Seidel methods, and others \cite{Hedy2013,Hedy2010}. To the best of our knowledge, the Kurdyka-{\L}ojasiewicz (KL) property plays a pivotal role in analyzing global sequential convergence. To continue, we review the following results.

\begin{lemma}\label{lema1}{\rm\cite{Hedy2013}}
Let $\varphi: \mathbb{R}^n \rightarrow \mathbb{R} \cup \{+ \infty\}$ be a proper lower semicontinuous function. Consider a sequence $\{\x^{(k)}\}_{k>0}\subseteq \mathbb{R}^n$ satisfying the following three conditions:

\noindent{\rm(i)}(Sufficient decrease condition) There exists $a>0$ such that
$$
\varphi(\x^{(k+1)})+a\|\x^{(k+1)}-\x^{(k)}\|^2\leq\varphi(\x^{(k)})
$$
holds for any $k\in\mathbb{N}$.

\noindent{\rm(ii)}(Relative error condition) There exist $b>0$ and $\omega^{(k+1)}\in\partial\varphi(\x^{(k+1)})$ (the subdifferential
 of $\varphi$ evaluated at $\x^{(k+1)}$), such that
$$
\|\omega^{(k+1)}\|\leq b\|\x^{(k+1)}-\x^{(k)}\|
$$
holds for any $k\in\mathbb{N}$.

\noindent{\rm(iii)} (Continuity condition) There exist subsequence $\{\x^{(k_j)}\}$ and $\x^*$ such that
$$
\x^{(k_j)}\rightarrow\x^*~~{\rm and}~~\varphi(\x^{(k_j)}) \rightarrow \varphi(\x^*),~~{\rm as}~j\rightarrow+\infty
$$
If $\varphi$ satisfies the KL property at $\x^*$, then  $\0\in\partial\varphi(\x^*)$, and $$\sum_{k=1}^{+\infty}\|\x^{(k+1)}-\x^{(k)}\|^2<+\infty,~~\lim_{k\rightarrow+\infty}\x^{(k)}=\x^*.$$
\end{lemma}

Denote the indicator function $\iota(\t)$ as below:
$$
\iota(\t)=\left\{
\begin{array}{ll}
0, &~~{\rm if}~\t\in \mathbb{D}_2,\\
+\infty, &~~{\rm otherwise}.\\
\end{array}
\right.
$$
Then the BMP \eqref{e9} can be equivalently reformulated as an unconstraint problem:
\begin{equation}\label{e15}
\min_{\t\in \mathbb{D}_2}F_{\alpha}(\t)=\min_{\t\in \mathbb{R}^m\times\mathbb{R}^n\times\mathbb{R}^m\times\mathbb{R}^n}F_{\alpha}(\t)+\iota(\t).
\end{equation}
Let $\varphi(\t)=F_{\alpha}(\t)+\iota(\t)$. Then, $\varphi(\t)$ holds the KL property automatically. From Theorem \ref{thm7}, we know that
the sufficient decrease condition and continuous condition are satisfied. Note that for all $\t\in \mathbb{D}_2$, it holds that
$$
\partial\varphi(\t)=\partial\left(F_{\alpha}(\t)+\iota(\t)\right)=\nabla F_{\alpha}(\t)+\partial\iota(\t),
$$
and $\partial\iota(\t)=\{\lambda \t : \lambda\in\mathbb{R}\}$, where $\partial$ means the Fr\'{e}chet subdifferential. Therefore, we know that $\omega^{(k+1)}=\nabla F_{\alpha}(\t^{(k+1)})-F_{\alpha}(\bar{\t})\t^{(k+1)}\in\partial\varphi(\t^{(k+1)})$ and the following result.

\begin{lemma}\label{lema2}
Let $\lambda^*=F_{\alpha}(\bar{\t})$, which is defined as in Theorem \ref{thm7}. Then, for all iterations $k\in\mathbb{N}$, there is an $M>0$ such that
$$
\|\omega^{(k+1)}\|=\|\nabla F_{\alpha}(\t^{(k+1)})-\lambda^*\t^{(k+1)}\|\leq M\|\t^{(k+1)}-\t^{(k)}\|.
$$
\end{lemma}
\proof
By Algorithm \ref{alg1} and the proof of Theorem \ref{thm7}, it's clear that
$$
\left\{
\begin{aligned}
&\|\nabla_{{\bf u}} F_{\alpha}({\bf u}^{(k+1)},\v^{(k)},\w^{(k)},\z^{(k)})-\lambda^*{\bf u}^{(k+1)}\|\leq \gamma_1\|{\bf u}^{(k+1)}-{\bf u}^{(k)}\|,\\
&\|\nabla_{\v}F_{\alpha}({\bf u}^{(k+1)},\v^{(k+1)},\w^{(k)},\z^{(k)})-\lambda^*\v^{(k+1)}\|\leq \gamma_2\|\v^{(k+1)}-\v^{(k)}\|,\\
&\|\nabla_{\w}F_{\alpha}({\bf u}^{(k+1)},\v^{(k+1)},\w^{(k+1)},\z^{(k)})-\lambda^*\w^{(k+1)}\|\leq \gamma_3\|\w^{(k+1)}-\w^{(k)}\|,\\
&\|\nabla_{\z}F_{\alpha}({\bf u}^{(k+1)},\v^{(k+1)},\w^{(k+1)},\z^{(k+1)})-\lambda^*\z^{(k+1)}\|\leq\gamma_4\|\z^{(k+1)}-\z^{(k)}\|,\\
\end{aligned}
\right.
$$
By the multilinear structure of $F_{\alpha}(\t)$ and the triangle inequality, we obtain that
$$
\begin{aligned}
&\|\nabla_{{\bf u}} F_{\alpha}(\t^{(k+1)})-\lambda^*{\bf u}^{(k+1)}\| \\
&\leq\|\nabla_{{\bf u}}  F_{\alpha}({\bf u}^{(k+1)},\v^{(k+1)},\w^{(k+1)}\z^{(k+1)})-\nabla_{{\bf u}}  F_{\alpha}({\bf u}^{(k+1)},\v^{(k+1)},\w^{(k+1)},\z^{(k)})\| \\
&+\|\nabla_{{\bf u}} F_{\alpha}({\bf u}^{(k+1)},\v^{(k+1)},\w^{(k+1)},\z^{(k)})-\nabla_{{\bf u}}  F_{\alpha}({\bf u}^{(k+1)},\v^{(k+1)},\w^{(k)},\z^{(k)})\|\\
&+\|\nabla_{{\bf u}} F_{\alpha}({\bf u}^{(k+1)},\v^{(k+1)},\w^{(k)},\z^{(k)})-\nabla_{{\bf u}}  F_{\alpha}({\bf u}^{(k+1)},\v^{(k)},\w^{(k)},\z^{(k)})\|\\
&+\|\nabla_{{\bf u}} F_{\alpha}({\bf u}^{(k+1)},\v^{(k)},\w^{(k)},\z^{(k)})-\nabla_{{\bf u}}  F_{\alpha}({\bf u}^{(k)},\v^{(k)},\w^{(k)},\z^{(k)})\|\\
&\leq m_1\|\t^{(k+1)}-\t^{(k)}\|,
\end{aligned}
$$
where $m_1\in\mathbb{R}$ is a constant. Similarly, there are positive constants $m_2, m_3, m_4,$ satisfying that, for each block of $\t^{(k+1)}$,
\begin{equation}\label{e16}
\left\{
\begin{aligned}
&\|\nabla_{{\bf u}} F_{\alpha}(\t^{(k+1)})-\lambda^*{\bf u}^{(k+1)}\|\leq m_1\|\t^{(k+1)}-\t^{(k)}\|, \\
&\|\nabla_{\v} F_{\alpha}(\t^{(k+1)})-\lambda^*\v^{(k+1)}\|\leq m_2\|\t^{(k+1)}-\t^{(k)}\|, \\
&\|\nabla_{\w} F_{\alpha}(\t^{(k+1)})-\lambda^*\w^{(k+1)}\|\leq m_3\|\t^{(k+1)}-\t^{(k)}\|, \\
&\|\nabla_{\z} F_{\alpha}(\t^{(k+1)})-\lambda^*\z^{(k+1)}\|\leq m_4\|\t^{(k+1)}-\t^{(k)}\|. \\
\end{aligned}
\right.
\end{equation}
Let $M=m_1+m_2+m_3+m_4$.  Combining \eqref{e16} with the definition of $\omega^{(k+1)}$, it follows that
$$
\begin{aligned}
\|\omega^{(k+1)}\|=&\|\nabla F_{\alpha}(\t^{(k+1)})-\lambda^*\t^{(k+1)}\|\\
\leq&\|\nabla_{{\bf u}} F_{\alpha}(\t^{(k+1)})-\lambda^*{\bf u}^{(k+1)}\|+\|\nabla_{\v} F_{\alpha}(\t^{(k+1)})-\lambda^*\v^{(k+1)}\| \\
&+\|\nabla_{\w} F_{\alpha}(\t^{(k+1)})-\lambda^*\w^{(k+1)}\|+\|\nabla_{\z} F_{\alpha}(\t^{(k+1)})-\lambda^*\z^{(k+1)}\|\\
\leq& M\|\t^{(k+1)}-\t^{(k)}\|,
\end{aligned}
$$
and the desired result follows.
\qed

By Theorem \ref{thm7}, Lemma \ref{lema1} and Lemma \ref{lema2}, we have the following global convergence.
\begin{theorem}\label{thm8}{\rm(Global sequence convergence)} Assume that the sequence $\{\t^{(k)}\}$ is an infinite sequence generated by Algorithm \ref{alg1}. Then, it converges globally to a KKT point $\bar{\t}$ of (\ref{e9}).
\end{theorem}

To end this section, we present a convergence rate of Algorithm \ref{alg1}.
First of all, we derive the KL exponent of the associated functions of $\varphi(\t)$.
To do this, we need the classical {\L}ojasiewicz gradient inequality for polynomials (See Theorem 4.2 of \cite{DK2005} for the detail).
\begin{lemma}\label{lema3}{\rm({\L}ojasiewicz gradient inequality)}
Let $h(\x)$ be a polynomial of $\mathbb{R}^n$ with degree $d\in\mathbb{N}$. Suppose that $h(\bar{\x})=0$. Then there exist constants $\epsilon, c>0$ such that, for all $\x\in\mathbb{R}^n$ with $\|\x-\bar{\x}\|\leq\epsilon$, we have
$$
\|\nabla h(\x)\|\geq c|h(\x)|^{1-\tau},~{\rm where}~\tau=R(n,d)^{-1},~R(n,d)=\left\{
\begin{array}{ll}
1,&{\rm if}~~ d=1,\\
d(3d-3)^{n-1},& {\rm if}~~ d\geq2.
\end{array}
\right.
$$
\end{lemma}

By Lemma \ref{lema3}, the following theorem establishes the explicit KL exponent of the merit function $\varphi(\t)$.
\begin{theorem}\label{thm9} Let $\varphi(\t)=F_{\alpha}(\t)+\iota(\t)$ be defined as in (\ref{e15}), and
$$\t=({\bf u}^\top, \v^\top, \w^\top, \z^\top)^\top\in\mathbb{D}_2\subseteq\R^{2(m+n)}.$$
Then, $\varphi(\t)$ satisfies the KL property with exponent $1-\tau$ at $\bar{\t}=(\bar{{\bf u}}^\top,\bar{\v}^\top,\bar{\w}^\top, \bar{\z}^\top)^\top$,
where $\tau=\frac{1}{4\times{81}^{m+n}}$ and $\bar{\t}$ is defined as in Theorem \ref{thm8}.
\end{theorem}
\proof
To prove the statement, let $\bar{\t}\in \mathbb{D}_2$, and let $\delta_1, \eta>0$ such that for any $\t$ satisfying $\|\t-\bar{\t}\|\leq \delta_1$, it follows that
$$
\varphi(\bar{\t})\leq \varphi(\t)\leq\varphi(\bar{\t})+\eta.
$$
On one hand, we can write
$$
\varphi(\t)=F_{\alpha}({\bf u},\v,\w,\t)+\iota({\bf u})+\iota(\v)+\iota(\w)+\iota(\z).
$$
By a direct computation, we know that
\begin{equation}\label{e17}
\left\{
\begin{aligned}
&\partial_{{\bf u}}\varphi(\t)=\{\nabla_{{\bf u}}F_{\alpha}(\t)+\lambda_1{\bf u} : \lambda_1\in\R\},\\
&\partial_{\v}\varphi(\t)=\{\nabla_{\v}F_{\alpha}(\t)+\lambda_2\v : \lambda_2\in\R\},\\
&\partial_{\w}\varphi(\t)=\{\nabla_{\w}F_{\alpha}(\t)+\lambda_3\v : \lambda_3\in\R\},\\
&\partial_{\z}\varphi(\t)=\{\nabla_{\z}F_{\alpha}(\t)+\lambda_4\z : \lambda_4\in\R\},\\
\end{aligned}
\right.
\end{equation}
which implies that
\begin{equation}\label{e18}
\begin{aligned}
&{\rm dist}(0, \partial\varphi(\t))^2=\inf_{\lambda_1,\lambda_2,\lambda_3,\lambda_4\in\R}\|\partial_{{\bf u}}\varphi(\t)\|^2+\|\partial_{\v}\varphi(\t)\|^2+\|\partial_{\w}\varphi(\t)\|^2+\|\partial_{\z}\varphi(\t)\|^2,\\
&=\inf_{\lambda_1\in\R}(\lambda_1^2+2\lambda_1F_{\alpha}(\t)+\|\nabla_{{\bf u}}F_{\alpha}(\t)\|^2)+\inf_{\lambda_2\in\R}(\lambda_2^2+2\lambda_2F_{\alpha}(\t)+\|\nabla_{\v}F_{\alpha}(\t)\|^2)\\
&+\inf_{\lambda_3\in\R}(\lambda_3^2+2\lambda_3F_{\alpha}(\t)+\|\nabla_{\w}F_{\alpha}(\t)\|^2)+\inf_{\lambda_4\in\R}(\lambda_4^2+2\lambda_4F_{\alpha}(\t)+\|\nabla_{\z}F_{\alpha}(\t)\|^2)\\
&=-4F_{\alpha}(\t)^2+\|\nabla F_{\alpha}(\t)\|^2.
\end{aligned}
\end{equation}

On the other hand, consider the following polynomial
$$
H({\bf u},\v,\w,\z,\lambda)=F_{\alpha}(\t)+\frac{\lambda}{2}(\|{\bf u}\|^2-1)+\frac{\lambda}{2}(\|\v\|^2-1)+\frac{\lambda}{2}(\|\w\|^2-1)
+\frac{\lambda}{2}(\|\z\|^2-1),$$
where $\lambda=-F_{\alpha}(\t)=-F_{\alpha}({\bf u},\v,\w,\z)$. Denote
$$\hat{H}({\bf u},\v,\w,\z,\lambda)=H({\bf u},\v,\w,\z,\lambda)-H(\bar{{\bf u}},\bar{\v},\bar{\w},\bar{\z},\bar{\lambda}),$$
where $\bar{\lambda}=-F_{\alpha}(\bar{{\bf u}},\bar{\v},\bar{\w},\bar{\z})$.
Obviously that $\hat{H}(\t)$ is a polynomial defined in $\R^{2(m+n)+1}$ with degree $4$ and $\hat{H}(\bar{{\bf u}},\bar{\v},\bar{\w},\bar{\z},\bar{\lambda})=0$.
By Lemma \ref{lema3}, there exist $\delta'>0, c>0$
such that, for all $\|\t-\bar{\t}\|\leq\delta_2$, it follows
$$
\begin{aligned}
\|\nabla H({\bf u},\v,\w,\z,\lambda)\|&=\|\nabla \hat{H}({\bf u},\v,\w,\z,\lambda)\|\\
&\geq c|H({\bf u},\v,\w,\z,\lambda)-H(\bar{{\bf u}},\bar{\v},\bar{\w},\bar{\z},\bar{\lambda})|^{1-\tau},
\end{aligned}
$$
where $\tau=\frac{1}{4\times 81^{m+n}}$.
Note that for any $\t\in\mathbb{D}_2, \lambda\in\R$, we have
$$
\left\{
\begin{aligned}
&\nabla_{{\bf u}} H({\bf u},\v,\w,\z,\lambda)=\nabla_{{\bf u}}F_{\alpha}(\t)+\lambda{\bf u},\\
&\nabla_{\v} H({\bf u},\v,\w,\z,\lambda)=\nabla_{\v}F_{\alpha}(\t)+\lambda\v,\\
&\nabla_{\w} H({\bf u},\v,\w,\z,\lambda)=\nabla_{\w}F_{\alpha}(\t)+\lambda\w,\\
&\nabla_{\z} H({\bf u},\v,\w,\z,\lambda)=\nabla_{\z}F_{\alpha}(\t)+\lambda\z,\\
&\nabla_{\lambda} H({\bf u},\v,\w,\z,\lambda)=0,\\
\end{aligned}
\right.
$$
which implies that
$$
\|\nabla H({\bf u},\v,\w,\z,\lambda)\|^2=-4F_{\alpha}(\t)^2+\|\nabla F_{\alpha}(\t)\|^2={\rm dist}(0, \partial\varphi(\t))^2,
$$
and $H({\bf u},\v,\w,\z,\lambda)=\varphi(\t)$, $H(\bar{{\bf u}},\bar{\v},\bar{\w},\bar{\z},\bar{\lambda})=\varphi(\bar{\t})$.
Take $\delta=\min\{\delta_1, \delta_2\}$. Combining this with (\ref{e18}), it holds hat, for all $\t\in\mathbb{D}_2$ with $\|\t-\bar{\t}\|\leq\delta$,
and $\varphi(\bar{\t})\leq \varphi(\t)\leq\varphi(\bar{\t})+\eta$,
$$
\begin{aligned}
{\rm dist}(0, \partial\varphi(\t))&\geq c|H({\bf u},\v,\w,\z,\lambda)-H(\bar{{\bf u}},\bar{\v},\bar{\w},\bar{\z},\bar{\lambda})|^{1-\tau}\\
&=c|\varphi(\t)-\varphi(\bar{\t})|^{1-\tau},
\end{aligned}
$$
and the desired results hold.
\qed
By Theorem \ref{thm9}, we now discuss the convergence rate of Algorithm \ref{alg1}. Under reasonable assumptions, we recall several classical results on convergence rates from the literature that established the convergence rate analyses based on the KL property \cite{Hedy2009,WCP18,YY2013}. Specifically, if the desingularization function of $\varphi(\t)$ is $\phi(s)=cs^{1-\alpha}$, then as shown in \cite{Hedy2009}, the following estimates hold:
\begin{itemize}
    \item[(i)] If $\alpha=0$, the sequence $\{\t^{(k)}\}$ converges in a finite number of steps.
    \item[(ii)] If $\alpha\in(0, 1/2]$, there exist constants $a>0$ and $\theta\in(0, 1)$ such that $\|\t^{(k)}-\bar{\t}\|\leq a\theta^k$.
    \item[(iii)] If $\alpha\in(1/2, 1)$, there exists a constant $a>0$ such that
    $$
    \|\t^{(k)}-\bar{\t}\|\leq ak^{-\frac{1-\alpha}{2\alpha-1}}.
    $$
\end{itemize}

Combining this with Theorem \ref{thm9} (where $\alpha=1-\tau$), we obtain the following result.

\begin{theorem}\label{thm10}
Assume $\{\t^{(k)}\}$ is an infinite sequence generated by Algorithm \ref{alg1}, and $\lim_{k\rightarrow\infty}\t^{(k)}=\bar{\t}$. Then, there exists a constant $a>0$ such that
$$
\|\t^{(k)}-\bar{\t}\|\leq ak^{-\tau/(1-2\tau)},
$$
where $\tau$ is defined in Theorem \ref{thm9}.
\end{theorem}

\proof
By Theorem \ref{thm9}, the merit function $\varphi(\t)$ satisfies the KL property with exponent $1-\tau$ at $\bar{\t}$. Hence, it suffices to verify that $1-\tau\in(1/2, 1)$. From the definition of $\tau$, it is evident that $0<\tau<1/2$, which ensures the desired result.
\qed

\section{Numerical results}\label{sec06}
In this section, we want to evaluate the practical performance of the proposed proximal alternating minimization (PAM) algorithm. Both biquadratic Cauchy tensor and general biquadratic tensor instances are considered in order to examine efficiency and numerical stability of PAM. All experiments are implemented in MATLAB R2018a and executed on a laptop equipped with a 13th Gen Intel(R) Core(TM) i7-13700H CPU (2.40 GHz) and 16 GB RAM.

To show the performance of PAM, we make comparison with maximum block improvement (MBI) method \cite{CSLZ2012}
and the alternating direction method (ADMM) \cite{Jiang}, which are representative block-type
algorithms widely used for polynomial optimization problems with spherical constraints.
In all experiments, the stopping criterion is defined as
\[
\mathrm{Err} :=
\frac{\left| f_{\alpha}(\x^{(k+1)},\y^{(k+1)}) - f_{\alpha}(\x^{(k)},\y^{(k)}) \right|}
{\max\left\{ \left| f_{\alpha}(\x^{(k+1)},\y^{(k+1)}) \right|,
\left| f_{\alpha}(\x^{(k)},\y^{(k)}) \right|, 1 \right\}}
\leq 10^{-6}.
\]
Unless otherwise stated, the proximal parameters in the PAM algorithm are set to
$\gamma_i = 0$ for $i = 1,2,3,4$ throughout all numerical experiments.
All algorithms are initialized with randomly generated feasible points that are uniformly
normalized onto the unit spheres.
The maximum number of iterations is set to $2000$ for all tested methods.

We first consider synthetic test problems generated from biquadratic Cauchy tensors in \eqref{e4}.
Specifically, the generating vectors
$\c, \d \in \mathbb{R}^n$ are randomly sampled with positive entries and then sorted in ascending
order.
The corresponding fourth-order tensor $\mathcal{A} = (a_{ijkl})$ is constructed as
\begin{equation}\notag
a_{ijkl} = \frac{1}{c_i + c_k + d_j + d_l}, \quad i,k\in [m], j,l \in [n].
\end{equation}
To ensure the concavity of the augmented objective function, we introduce a regularization term with parameter
$\alpha$ in accordance with the sufficient concavity condition derived from the theoretical analysis. The resulting tensor is then symmetrized to align with the multilinear reformulation. Our experiments cover a range of dimensions
$
n \in \{ 5,10,20,30,40,50,60,70,80,100 \}.
$
For each dimension, we perform 10 independent trials using randomly generated initial points. Due to the randomness of the generated test instances, which may not always meet the theoretical convergence criteria, we report the success ratio (labeled ``Succ."), defined as the proportion of trials that converge successfully within the prescribed maximum number of iterations.

\begin{table}[htbp]
\centering
\caption{Performance of PAM, MBI, and ADMM on randomly generated biquadratic Cauchy tensor instances}
\label{tab:rand_case}
\resizebox{\textwidth}{!}{
\fontsize{16pt}{16pt}\selectfont
\setlength{\tabcolsep}{7pt}
\renewcommand{\arraystretch}{1.6}
\begin{tabular}{c cccc cccc cccc}
\toprule
\bfseries
 & \multicolumn{4}{c}{PAM}
 & \multicolumn{4}{c}{MBI}
 & \multicolumn{4}{c}{ADMM} \\
\cmidrule(lr){2-5} \cmidrule(lr){6-9} \cmidrule(lr){10-13}
\bfseries
$n$
& Iter & Time & Obj. & Succ.
& Iter & Time & Obj. & Succ.
& Iter & Time & Obj. & Succ. \\
\midrule
\normalfont
$5$   & 9.0 & 0.04 & $-16.13$   & 1.00 & 21.8 & 0.03 & $-16.13$   & 1.00 & 57.0 & 0.09 & $-16.44$   & 0.90 \\
$10$  & 5.8 & 0.02 & $-60.58$   & 1.00 & 14.4 & 0.05 & $-60.58$   & 1.00 & 45.1 & 0.16 & $-60.58$   & 1.00 \\
$20$  & 4.2 & 0.02 & $-213.48$  & 1.00 & 11.7 & 0.05 & $-213.48$  & 1.00 & 81.2 & 0.40 & $-213.48$  & 1.00 \\
$30$  & 4.0 & 0.04 & $-533.81$  & 1.00 & 13.5 & 0.09 & $-476.68$  & 1.00 & 27.7 & 0.21 & $-534.42$  & 0.90 \\
$40$  & 4.5 & 0.18 & $-987.80$  & 1.00 & 13.8 & 0.37 & $-987.80$  & 1.00 & 22.6 & 0.57 & $-969.17$  & 0.70 \\
$50$  & 4.4 & 0.40 & $-1595.08$ & 1.00 & 13.7 & 0.94 & $-1595.07$ & 1.00 & 20.3 & 1.40 & $-1602.95$ & 0.90 \\
$60$  & 3.8 & 0.61 & $-2142.63$ & 1.00 & 12.2 & 1.59 & $-2142.63$ & 1.00 & 20.4 & 2.70 & $-2125.67$ & 0.80 \\
$70$  & 4.3 & 0.79 & $-2974.21$ & 1.00 & 12.5 & 1.90 & $-2974.21$ & 1.00 & 24.0 & 3.81 & $-2970.18$ & 0.70 \\
$80$  & 4.0 & 1.01 & $-3960.43$ & 1.00 & 11.9 & 2.44 & $-3960.43$ & 1.00 & 21.6 & 4.66 & $-3991.98$ & 0.70 \\
$100$ & 4.2 & 2.21 & $-6091.48$ & 1.00 & 12.2 & 4.96 & $-6091.48$ & 1.00 & 23.5 & 8.23 & $-6140.89$ & 0.60 \\
\bottomrule
\end{tabular}}
\end{table}

Table 1 presents the averaged numerical results of PAM, MBI, and ADMM with respect to the number of iterations (Iter), CPU time in seconds (Time), objective value obtained (Obj.), and success ratio (Succ.). It can be seen that PAM consistently requires significantly fewer iterations than both MBI and ADMM across all dimensions tested. Notably, the iteration count of PAM remains nearly stable as the problem dimension increases, while MBI and ADMM show a clear rise in the number of iterations. In terms of computational efficiency, PAM is substantially faster than the other two methods on medium- and large-scale instances, highlighting its advantage in exploiting the multi-block structure of the multilinear reformulation.
Moreover, PAM attains a success rate of $100\%$ in all cases tested, demonstrating robust convergence behavior for biquadratic Cauchy tensor instances.

To further demonstrate that the performance of the proposed PAM algorithm is not limited to specially structured tensor instances, we further investigate its numerical behavior on general BPPs. The corresponding results are summarized in Table 2. In these experiments, the only difference from the previous setting lies in the construction of tensor $\A$, while all other parameters and stopping criteria remain consistent with those used in Table 1.

\begin{table}[htbp]
\centering
\caption{Performance of PAM, MBI, and ADMM on General Random Biquadratic Tensor Instances}
\label{tab:rand_case}
\resizebox{\textwidth}{!}{
\fontsize{12pt}{14pt}\selectfont
\setlength{\tabcolsep}{6pt}
\renewcommand{\arraystretch}{1.4}
\begin{tabular}{c cccc cccc cccc}
\toprule
 & \multicolumn{4}{c}{PAM}
 & \multicolumn{4}{c}{MBI}
 & \multicolumn{4}{c}{ADMM} \\
\cmidrule(lr){2-5} \cmidrule(lr){6-9} \cmidrule(lr){10-13}
$n$
& Iter & Time & Obj. & Succ.
& Iter & Time & Obj. & Succ.
& Iter & Time & Obj. & Succ. \\
\midrule
$5$   & 11.0 & 0.02 & $-11.36$   & 1.00 & 56.1 & 0.07 & $-12.55$   & 1.00 & 307.6 & 0.49 & $-12.61$   & 0.50 \\
$10$  &  7.2 & 0.02 & $-49.99$   & 1.00 & 16.6 & 0.04 & $-49.99$   & 1.00 &  19.0 & 0.04 & $-49.66$   & 0.30 \\
$20$  & 10.4 & 0.05 & $-200.09$  & 1.00 & 12.5 & 0.05 & $-200.09$  & 1.00 &  20.1 & 0.09 & $-200.16$  & 0.70 \\
$30$  &  4.6 & 0.04 & $-449.94$  & 1.00 &  9.9 & 0.07 & $-449.94$  & 1.00 &  23.4 & 0.18 & $-449.78$  & 0.50 \\
$40$  & 11.6 & 0.20 & $-800.15$  & 1.00 & 13.6 & 0.24 & $-800.15$  & 1.00 &  21.0 & 0.37 & $-800.07$  & 0.50 \\
$50$  &  4.0 & 0.20 & $-1249.80$ & 1.00 &  9.0 & 0.39 & $-1249.80$ & 1.00 &  21.3 & 0.76 & $-1249.75$ & 0.40 \\
$60$  &  3.7 & 0.37 & $-1799.75$ & 1.00 & 13.3 & 1.20 & $-1799.75$ & 1.00 &  20.7 & 2.13 & $-1799.65$ & 0.30 \\
$70$  &  3.2 & 0.45 & $-2450.00$ & 1.00 &  9.3 & 1.05 & $-2450.00$ & 1.00 &  20.3 & 1.87 & $-2449.85$ & 0.30 \\
$80$  &  4.9 & 1.52 & $-3200.00$ & 1.00 & 11.0 & 3.33 & $-3200.00$ & 1.00 &  20.6 & 5.80 & $-3199.98$ & 0.70 \\
$100$ &  3.2 & 2.58 & $-4999.95$ & 1.00 &  9.2 & 5.81 & $-4999.95$ & 1.00 &  22.3 & 12.20& $-4999.88$ & 0.40 \\
\bottomrule
\end{tabular}}
\end{table}

Fig. 1 compares the numerical performance of the proposed PAM algorithm with the MBI method on biquadratic tensor instances of increasing dimensions, ranging from
$n=5$ to $n=100$. The first set of subfigures presents results for biquadratic Cauchy tensor instances, whereas the second set corresponds to general biquadratic tensor instances. In both cases, the regularized problem  $F_\alpha$ is solved with a fixed parameter $\alpha = 6$, and the reported results are averaged across multiple independent trials for each dimension.

\begin{figure}[htbp]
  \centering

  % ËõÐ¡ figure ÓëÕýÎÄ/Ò³ÃæµÄ´¹Ö±¼ä¾à£¨¾Ö²¿ÉúÐ§£©
  \setlength{\intextsep}{2pt}
  \setlength{\textfloatsep}{2pt}

  % ËõÐ¡×ÓÍ¼ caption ÓëÍ¼Æ¬Ö®¼äµÄ¾àÀë
  \setlength{\abovecaptionskip}{1pt}
  \setlength{\belowcaptionskip}{09pt}

  % -------- first row --------
  \begin{subfigure}[t]{0.48\textwidth}
    \centering
    \includegraphics[width=1.08\linewidth, trim=1cm 10cm 2cm 10cm, clip]{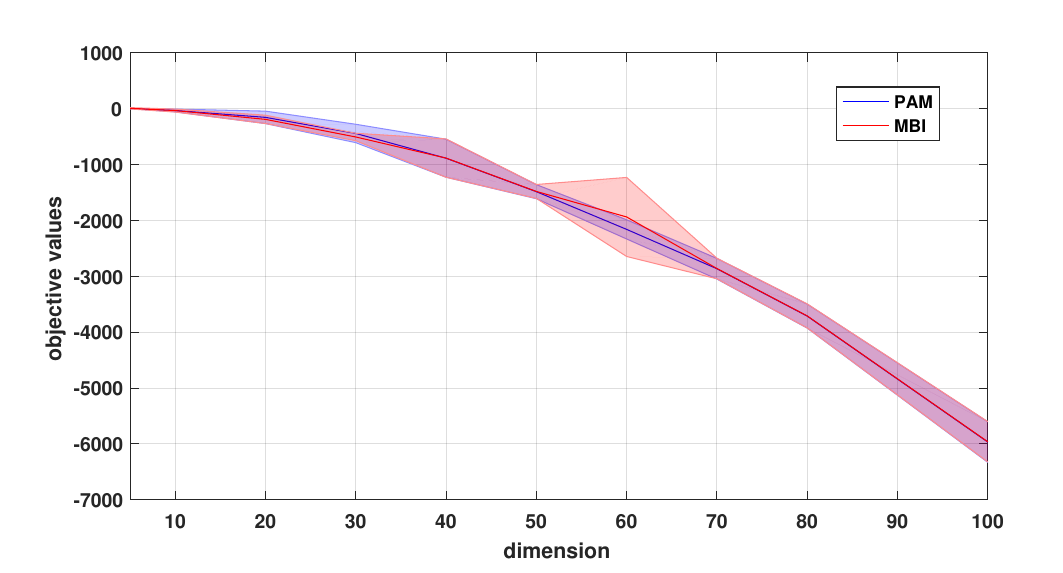}
    \caption{}
    \label{fig:fig1a}
  \end{subfigure}
  \hspace{0.02\textwidth}
  \begin{subfigure}[t]{0.48\textwidth}
    \centering
    \includegraphics[width=1.08\linewidth, trim=1cm 10cm 2cm 10cm, clip]{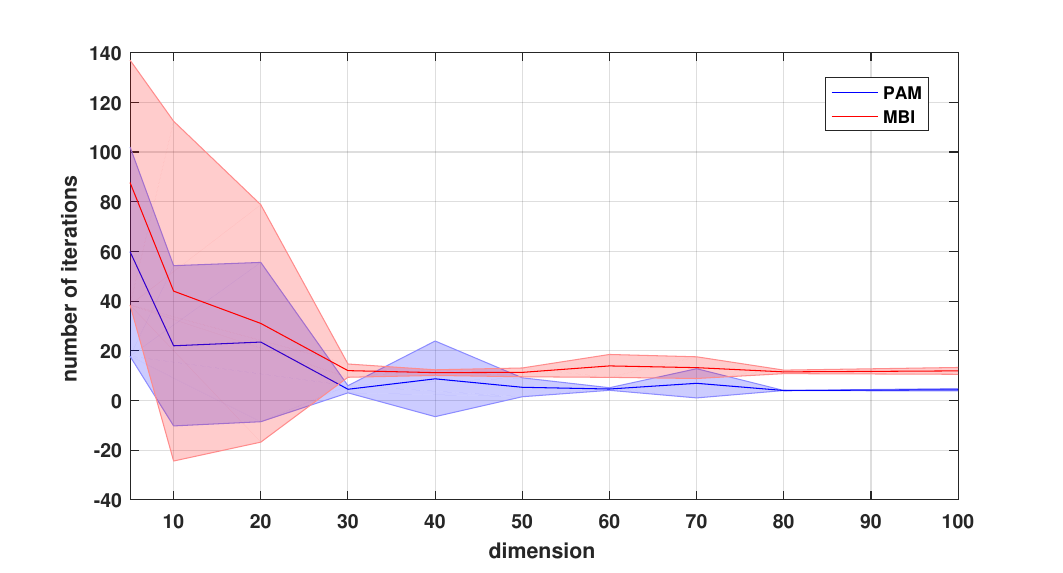}
    \caption{}
    \label{fig:fig1b}
  \end{subfigure}

  \vspace{-0.4cm}  % ¿ØÖÆÁ½ÐÐ×ÓÍ¼Ö®¼äµÄ¾àÀë£¨²»ÒªÔÙÓÃ -0.9cm£©

  % -------- second row --------
  \begin{subfigure}[t]{0.48\textwidth}
    \centering
    \includegraphics[width=1.08\linewidth, trim=1cm 10cm 2cm 10cm, clip]{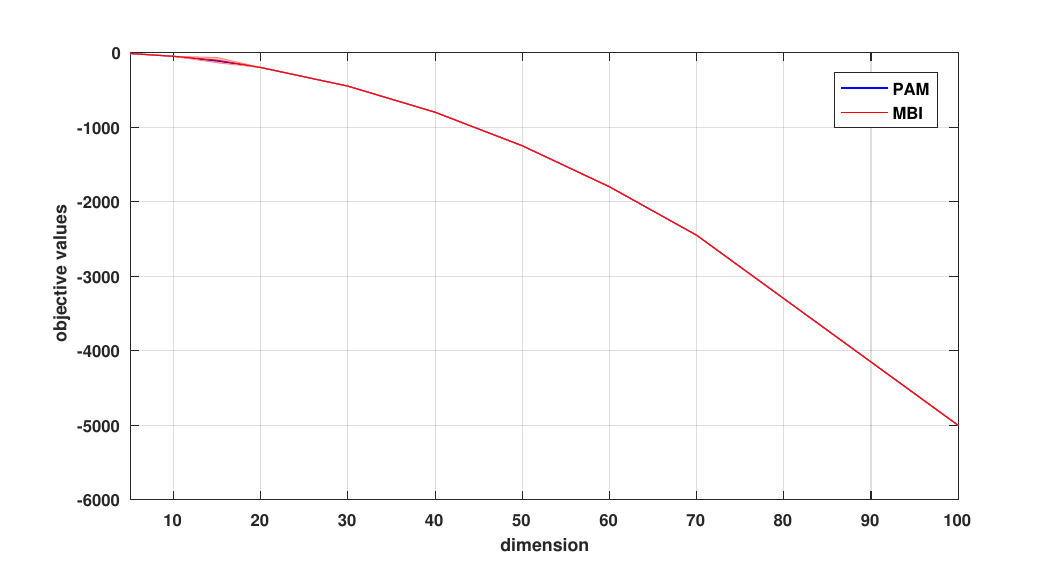}
    \caption{}
    \label{fig:fig1c}
  \end{subfigure}
  \hspace{0.02\textwidth}
  \begin{subfigure}[t]{0.48\textwidth}
    \centering
    \includegraphics[width=1.08\linewidth, trim=1cm 10cm 2cm 10cm, clip]{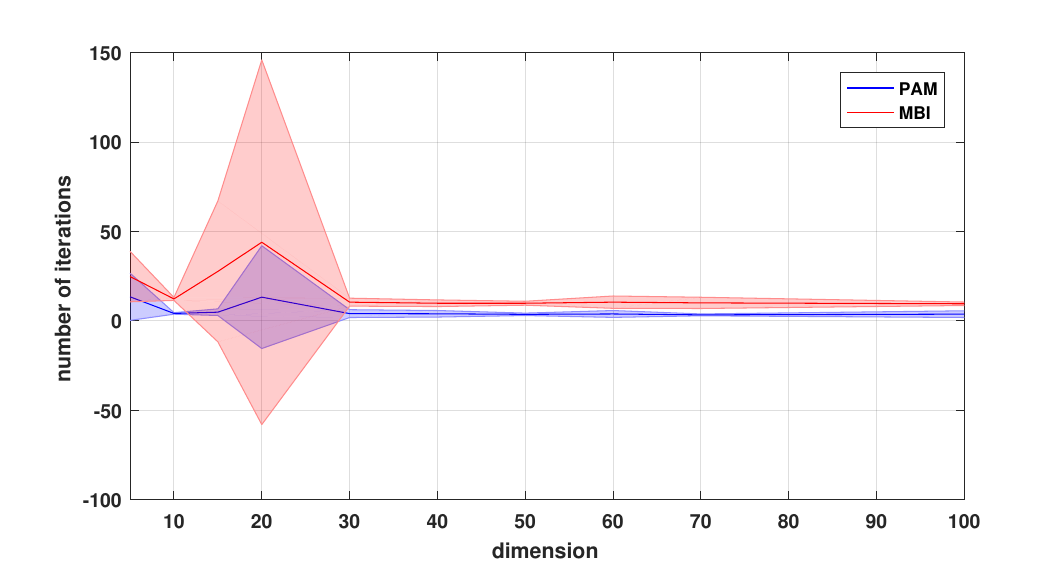}
    \caption{}
    \label{fig:fig1d}
  \end{subfigure}

  \caption{Performance comparison between PAM and MBI across different dimensions.
  Subfigures (a),(b) correspond to biquadratic Cauchy tensor instances, while subfigures (c),(d) report the results obtained on general biquadratic tensor instances.}
  \label{fig:fig1}
\end{figure}

In Fig.1(a) and (c), the average objective values consistently decrease as the problem dimension grows. In both cases, PAM achieves objective values comparable to¡ªor slightly better than¡ªthose of MBI, while exhibiting a narrower variation range, particularly at higher dimensions. Fig. 1(b) and (d) report the average number of iterations needed for convergence. For biquadratic Cauchy tensor instances, PAM requires fewer iterations than MBI at low and moderate dimensions and remains stable as dimension increases, whereas MBI shows larger fluctuations. A similar pattern is observed for general biquadratic tensor instances: PAM converges in a low, nearly dimension-invariant number of iterations, while MBI demands more iterations and displays greater variability. Overall, these results demonstrate that PAM is more efficient and numerically stable than MBI¡ªan advantage that holds for both structured and general biquadratic tensor problems.

Fig. 2 illustrates the convergence behavior of the proposed PAM algorithm under different initialization strategies. We consider two types of biquadratic tensor instances: biquadratic Cauchy tensors and general biquadratic tensors. In both cases, the tensor $\mathcal{A}$ is constructed accordingly and held fixed, with the problem dimension set to $n=100$ and the regularization parameter fixed at $\alpha=6$. Initial points for $(u,v,w,z)$ are drawn from Gaussian and uniform distributions, normalized onto their respective unit spheres, and tested over 20 independent runs.

\begin{figure}[htbp]
  \centering
  \begin{subfigure}[t]{0.48\textwidth}
    \centering
    \includegraphics[width=1.15\linewidth, trim=1cm 10cm 2cm 10cm, clip]{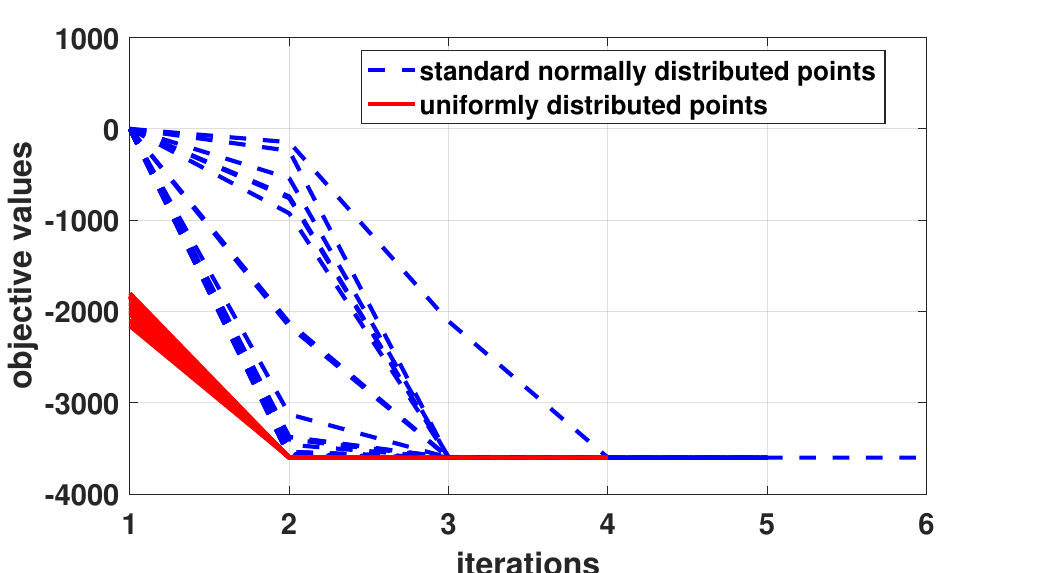}
    \caption{}
    \label{fig:fig2a}
  \end{subfigure}
  \hfill
  \begin{subfigure}[t]{0.49\textwidth}
    \centering
    \includegraphics[width=1.15\linewidth, trim=1cm 10cm 2cm 10cm, clip]{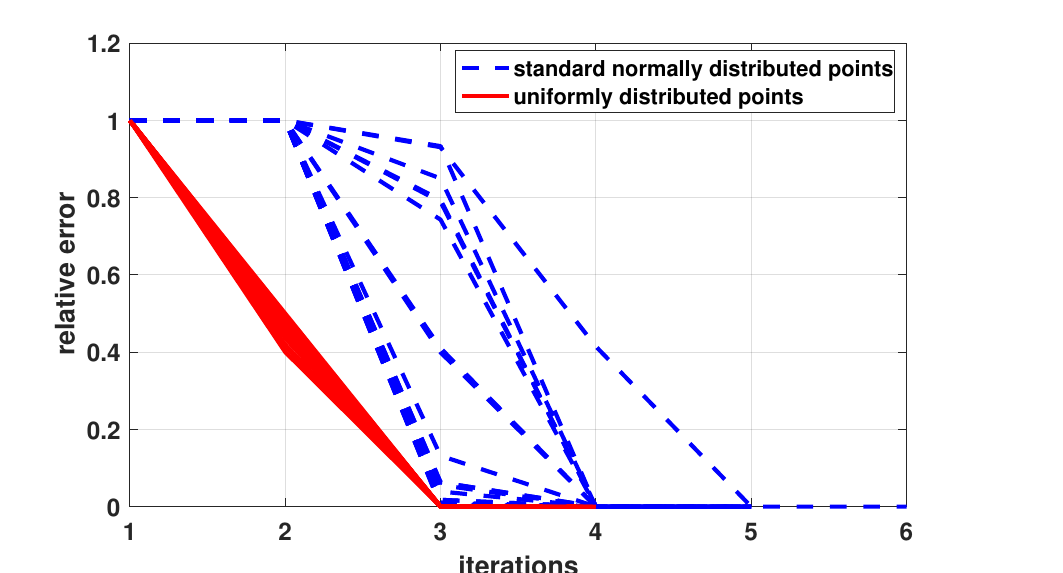}
    \caption{ }
    \label{fig:fig2b}
  \end{subfigure}
    \begin{subfigure}[t]{0.48\textwidth}
    \centering
    \includegraphics[width=1.15\linewidth, trim=1cm 10cm 2cm 10cm, clip]{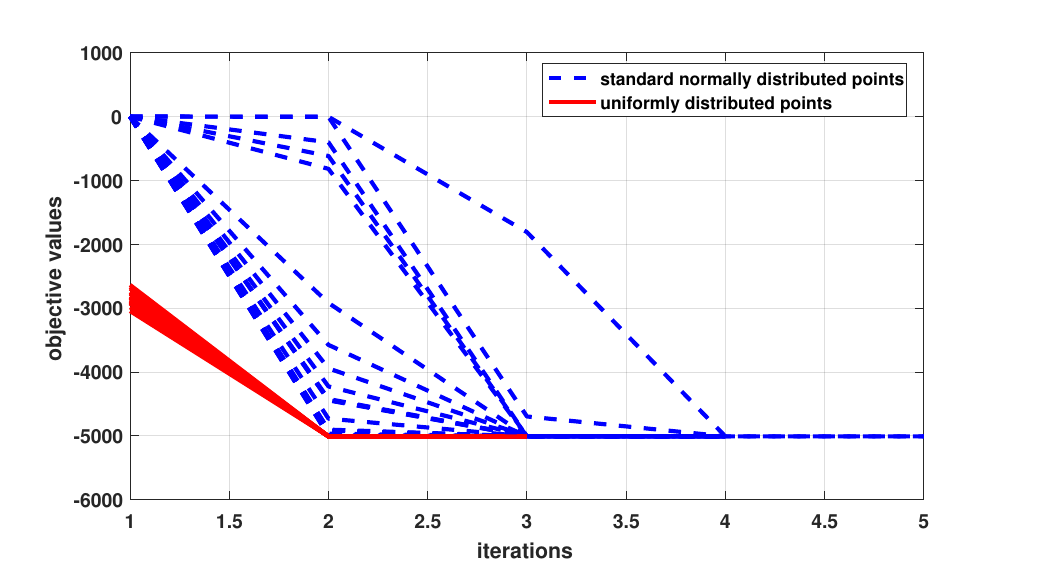}
    \caption{}
    \label{fig:fig2a}
  \end{subfigure}
  \hfill
  \begin{subfigure}[t]{0.49\textwidth}
    \centering
    \includegraphics[width=1.15\linewidth, trim=1cm 10cm 2cm 10cm, clip]{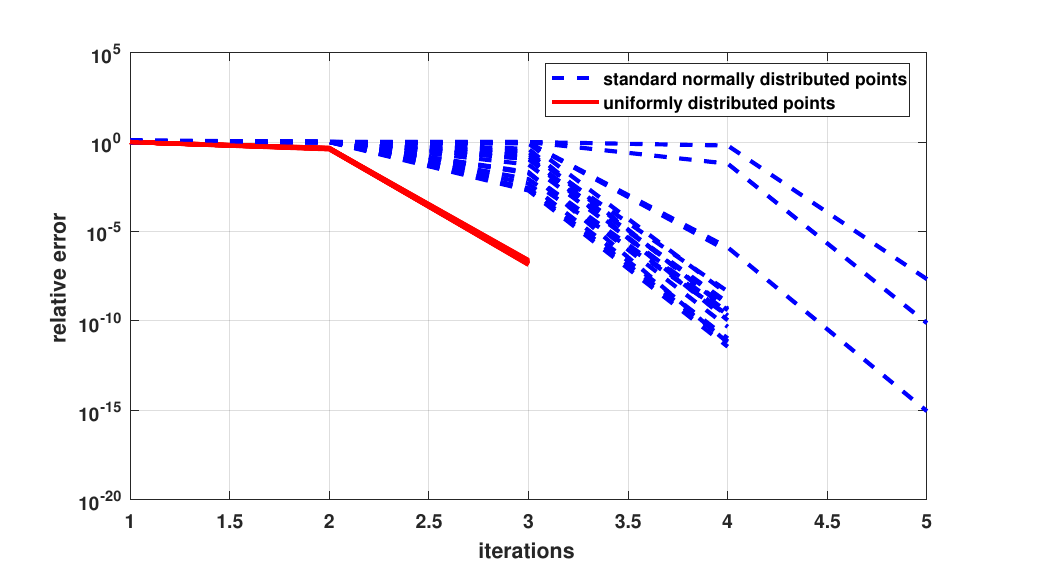}
    \caption{ }
    \label{fig:fig2b}
  \end{subfigure}
  \caption{Convergence behavior of the PAM algorithm under different initialization strategies, showing the evolution of both objective values and relative stopping errors across iterations}
  \label{fig:fig2}
\end{figure}

As shown in Fig. 2(a) and (c), the objective values decrease rapidly in early iterations before stabilizing, a trend consistent across both initialization strategies and for both Cauchy and general biquadratic tensor instances. Correspondingly, Fig. 2(b) and (d) show that the relative stopping errors exhibit an overall decaying trend, with occasional mild non-monotonicity due to the blockwise update scheme of PAM. Compared to Gaussian initialization, uniformly distributed starting points yield slightly more concentrated convergence curves, indicating marginally better numerical stability.

Recognizing that the multilinear reformulation employs an augmented (regularized) objective function, we want to further examine the sensitivity of the proposed PAM algorithm to the regularization parameter $\alpha$. The experiment considers two types of biquadratic tensor instances: biquadratic Cauchy tensors and general biquadratic tensors. For both cases, we test values of $\alpha \in \{1,5,15,25,35\}$, while keeping all other algorithmic settings fixed, and vary the problem dimension from $n=5$ to $n=100$.

As illustrated in Fig. 3(a) and (c), the average objective values attained under different $\alpha$ settings remain closely aligned across both Cauchy and general biquadratic tensor instances, consistently decreasing as the problem dimension grows. This indicates that the solution quality of the proposed PAM algorithm is largely insensitive to the choice of $\alpha$. Correspondingly, Fig. 3(b) and (d) present the associated CPU times. While larger $\alpha$ values may moderately raise computational cost at higher dimensions, the overall growth pattern remains similar across all parameter settings. Together, these observations underscore the robustness of PAM with respect to the regularization parameter for both structured and general biquadratic tensor problems.

\begin{figure}[htbp]
  \centering
  \begin{subfigure}[t]{0.49\textwidth}
    \centering
    \includegraphics[width=1.1\linewidth, trim=1cm 10cm 2cm 10cm, clip]{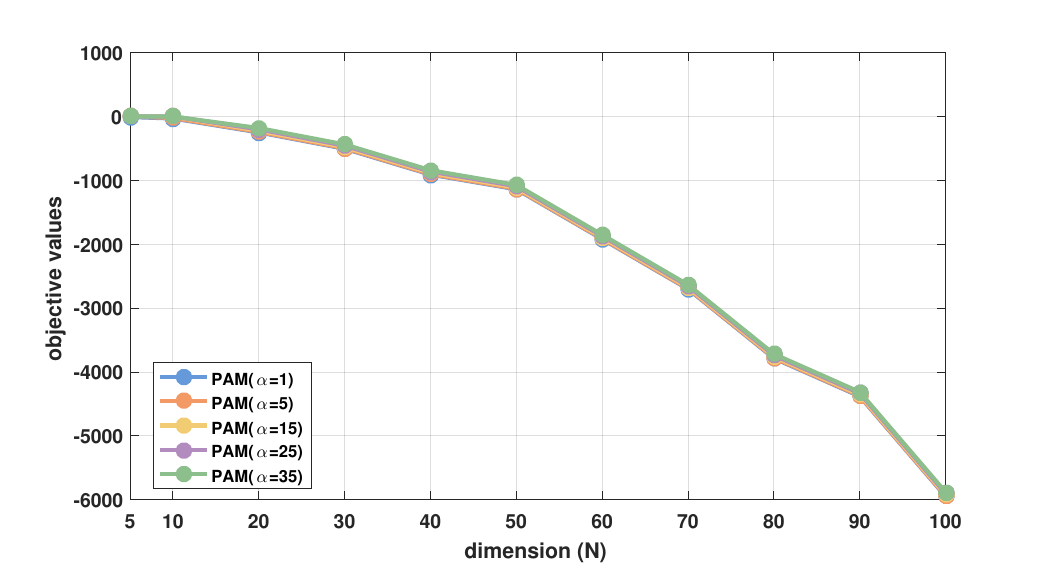}
    \caption{}
    \label{fig:fig3a}
  \end{subfigure}
  \hfill
  \begin{subfigure}[t]{0.49\textwidth}
    \centering
    \includegraphics[width=1.1\linewidth, trim=1cm 10cm 2cm 10cm, clip]{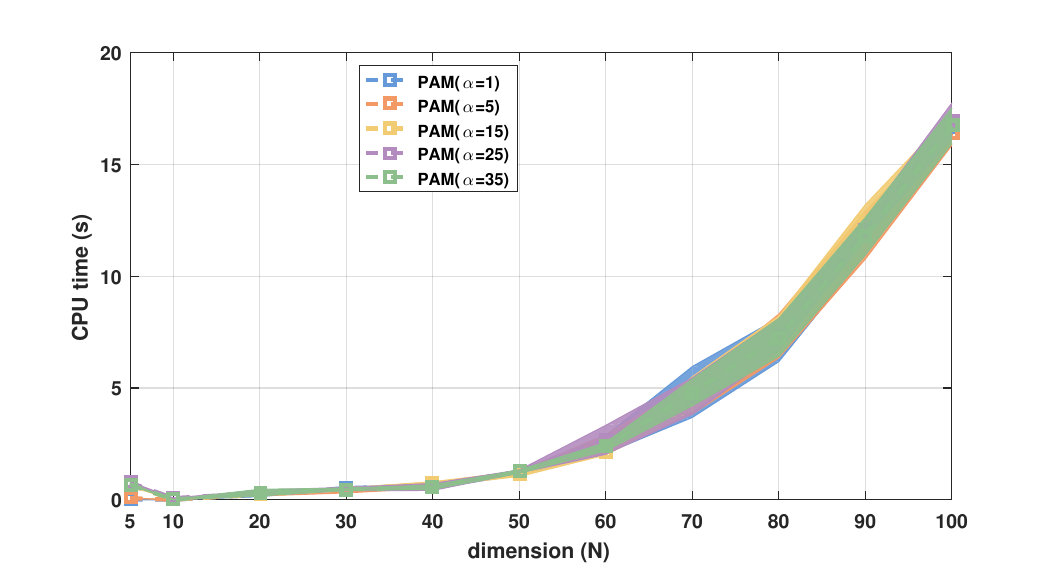}
    \caption{ }
    \label{fig:fig3b}
  \end{subfigure}
  \begin{subfigure}[t]{0.49\textwidth}
    \centering
    \includegraphics[width=1.1\linewidth, trim=1cm 10cm 2cm 10cm, clip]{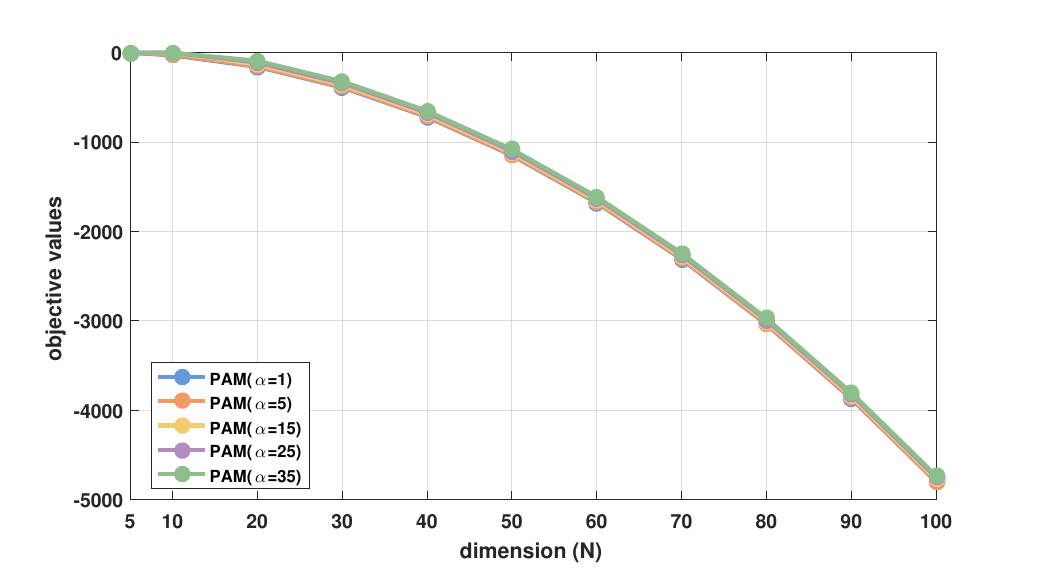}
    \caption{}
    \label{fig:fig3a}
  \end{subfigure}
  \hfill
  \begin{subfigure}[t]{0.49\textwidth}
    \centering
    \includegraphics[width=1.1\linewidth, trim=1cm 10cm 2cm 10cm, clip]{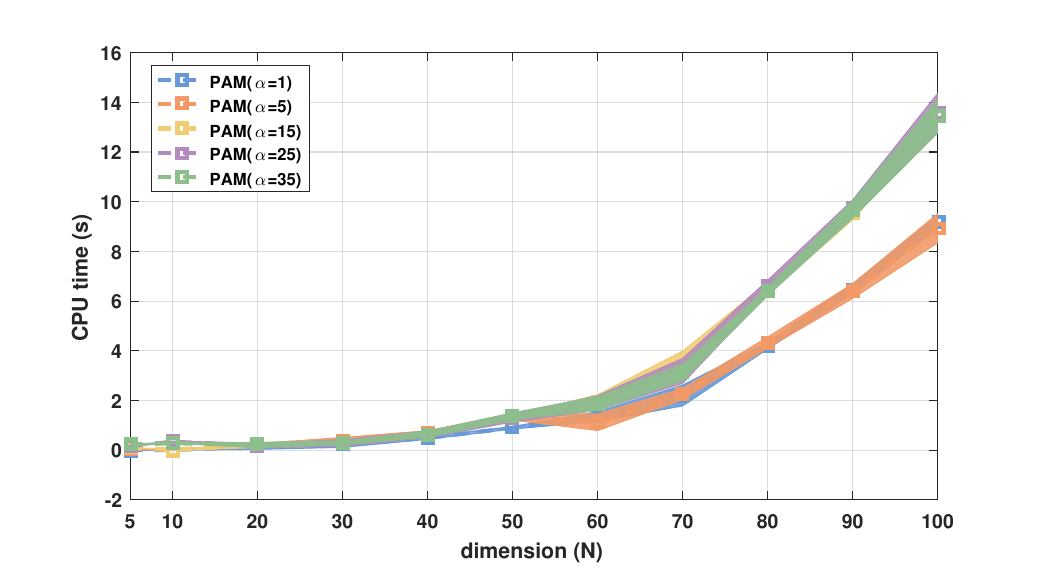}
    \caption{ }
    \label{fig:fig3b}
  \end{subfigure}
  \caption{Sensitivity Analysis of the PAM Algorithm with Respect to $\alpha$}
  \label{fig:fig3}
\end{figure}

Finally, since the proposed PAM algorithm involves a proximal parameter $\gamma$ that influences both convergence behavior and computational efficiency, we conduct a detailed sensitivity analysis. Two classes of biquadratic tensor instances are examined: biquadratic Cauchy tensors and general biquadratic tensors. In both cases, the regularization parameter is fixed at  $\alpha = 6$ while $\gamma$ varies over $
\gamma \in \{0,\,0.2,\,0.4,\,0.6,\,0.8,\,1\}$, with all other algorithmic settings held constant. The problem dimension is gradually increased from
$n=5$ to $n=100$.
Figs. 4(e) and 4(f) reveal that the obtained objective values remain similar across different values of $\gamma$, suggesting limited impact on solution quality. In contrast, Figs. 4(a) and 4(d) show that iteration counts and CPU times depend more noticeably on $\gamma$ for small to moderate dimensions. A larger $\gamma$ generally improves the likelihood of converging to the minimum eigenvalue, albeit at increased computational expense. As the dimension grows, these differences diminish, and the performance curves for different
$\gamma$ become nearly indistinguishable.

Overall, the proposed PAM algorithm exhibits considerable robustness to the choice of $\gamma$. This parameter primarily affects the convergence behavior rather than the final solution quality. While a smaller $\gamma$ often leads to faster convergence, it increases the risk of stopping at local solutions. Conversely, a larger $\gamma$ trades off higher computational cost for a greater likelihood of converging to the minimum eigenvalue. As this trade-off diminishes in higher-dimensional problems, we set $\gamma = 0$ for all experiments in this paper to ensure stable performance and avoid extra tuning.

\begin{figure}[htbp]
  \centering
  \begin{subfigure}[t]{0.49\textwidth}
    \centering
    \includegraphics[width=1.1\linewidth, trim=1cm 10cm 2cm 10cm, clip]{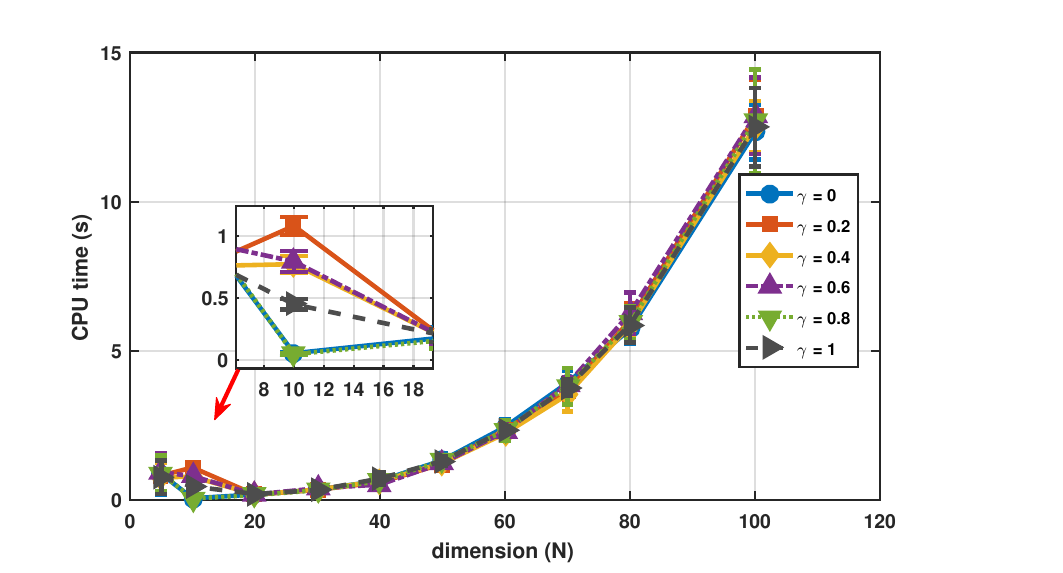}
    \caption{}
    \label{fig:fig3a}
  \end{subfigure}
  \hfill
  \begin{subfigure}[t]{0.49\textwidth}
    \centering
    \includegraphics[width=1.1\linewidth, trim=1cm 10cm 2cm 10cm, clip]{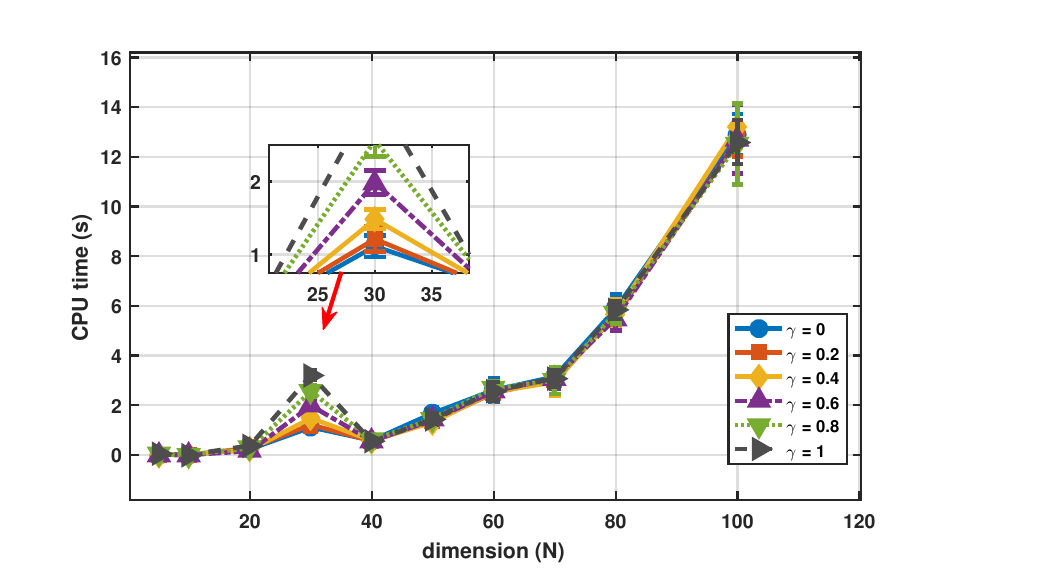}
    \caption{ }
    \label{fig:fig3b}
  \end{subfigure}
  \begin{subfigure}[t]{0.49\textwidth}
    \centering
    \includegraphics[width=1.1\linewidth, trim=1cm 10cm 2cm 10cm, clip]{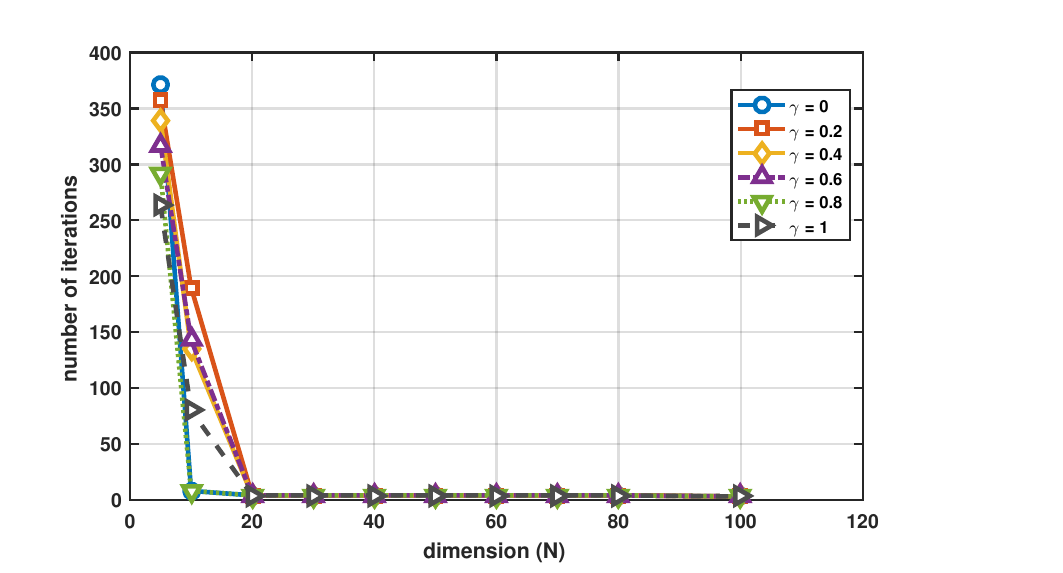}
    \caption{}
    \label{fig:fig3a}
  \end{subfigure}
  \hfill
  \begin{subfigure}[t]{0.49\textwidth}
    \centering
    \includegraphics[width=1.1\linewidth, trim=1cm 10cm 2cm 10cm, clip]{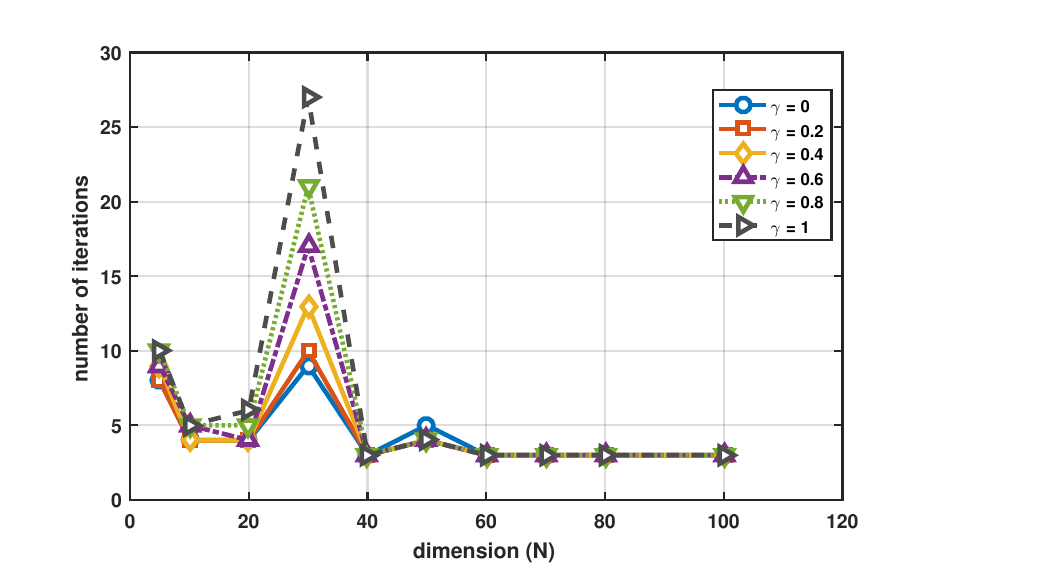}
    \caption{ }
    \label{fig:fig3b}
  \end{subfigure}
  \begin{subfigure}[t]{0.49\textwidth}
    \centering
    \includegraphics[width=1.1\linewidth, trim=1cm 10cm 2cm 10cm, clip]{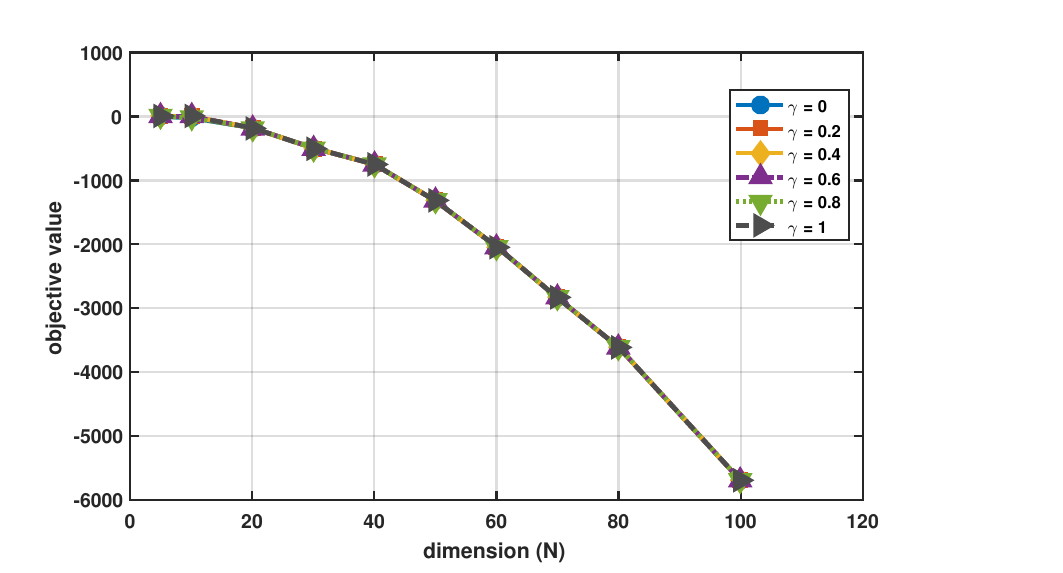}
    \caption{}
    \label{fig:fig3a}
  \end{subfigure}
  \hfill
  \begin{subfigure}[t]{0.49\textwidth}
    \centering
    \includegraphics[width=1.1\linewidth, trim=1cm 10cm 2cm 10cm, clip]{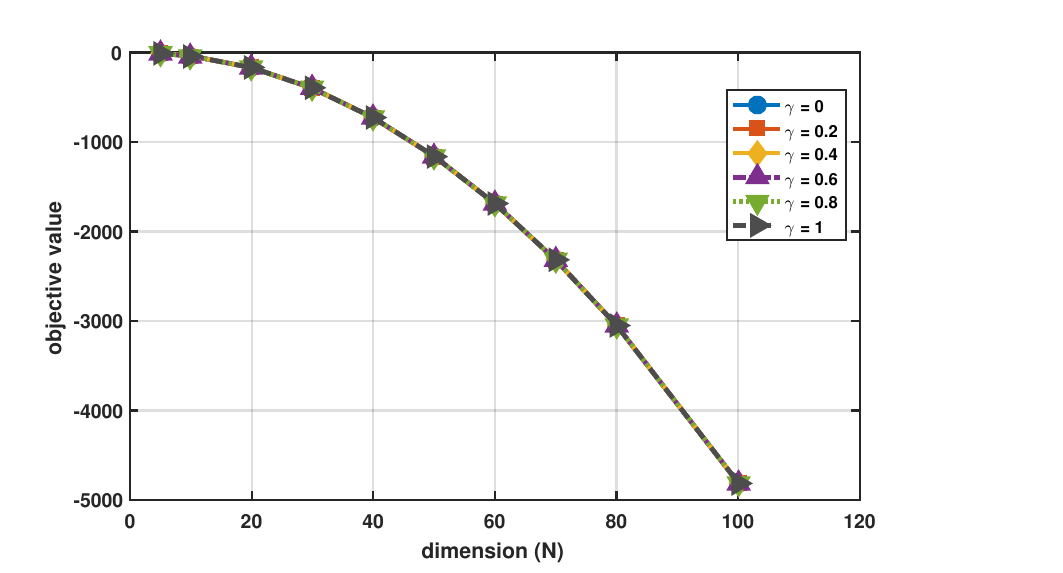}
    \caption{ }
    \label{fig:fig3b}
  \end{subfigure}
  \caption{Sensitivity of the PAM algorithm to the proximal parameter $\gamma$: CPU time, iteration count, and objective value for biquadratic Cauchy (left) and general biquadratic (right) tensor instances.}
  \label{fig000}
\end{figure}

\section{Conclusions}\label{sec_con}

In this paper, we first establish several necessary and sufficient conditions for the positive semi-definiteness and positive definiteness of biquadratic Cauchy tensors. Leveraging the structured properties of these tensors, we then prove that the BPP and its equivalent multilinear formulation share the same set of optimal solutions. Then, the global sequence convergence of the PAM algorithm via the Kurdyka-{\L}ojasiewicz (KL) property is established. Furthermore, by reformulating the equivalent multilinear problem as an unconstrained optimization model, we enable the analysis of its KL exponent and derive an explicit expression for the convergence rate of PAM. Finally, numerical experiments are conducted on both biquadratic Cauchy tensors and general biquadratic tensor instances to evaluate the efficiency, stability, and practical performance of the proposed algorithm.

Several intriguing and open questions remain for future investigation. In practical applications, we always face spherical optimization problems with inhomogeneous polynomials. So, does a similar equivalence with a multilinear optimization model still hold? For the homogeneous case, can the spherical constraint be generalized to other compact constraint sets? On the other hand, Combining tensor decomposition techniques may improve computational efficiency for certain practical problems. In the future, leveraging the low-rank and sparse structures of tensors, we will focus on
discussing algorithm acceleration as well as analyzing the computational complexity of improved algorithms with the help of various tensor decomposition techniques such as TT decomposition \cite{AG2024,IVO2011}, Tucker decomposition and CP decomposition \cite{Kolda2009}. Although these questions present significant theoretical challenges,
their resolution would yield valuable insights for implementable optimization methods.

\medskip
\begin{acknowledgements}
%The authors gratefully acknowledge the editor and two anonymous referees for their constructive comments, which significantly contributed to improving the quality of this paper.
This work was partially supported by the National Natural Science Foundation of China (12071249), Shandong Provincial Natural Science Foundation (ZR2024MA 003,ZR2021JQ01), and
research Center for Intelligent Operations Research, The Hong Kong Polytechnic University (4-ZZT8),
\end{acknowledgements}

%\noindent{\bf Availability of data and materials:} No datasets were generated or analysed during the current study.

\end{document}